\documentclass[12pt]{amsart}

\usepackage{epsf,amssymb,stmaryrd}

\newtheorem{thm}{Theorem}[section]

\newtheorem{prop}[thm]{Proposition}
\newtheorem{defn}[thm]{Definition}
\newtheorem{lemma}[thm]{Lemma}

\newtheorem*{op}{Operation}

\newcommand{\R}{\mathbb{R}}
\newcommand{\Z}{\mathbb{Z}}

\newcommand{\bdry}{\partial}

\newcommand{\s}{\vskip.1in}
\newcommand{\n}{\noindent}
\newcommand{\ii}{\mbox{\it II}}
\renewcommand{\ll}{\llbracket}
\newcommand{\rr}{\rrbracket}

\newcommand{\be}{\begin{enumerate}}
\newcommand{\ee}{\end{enumerate}}

\begin{document}
\title{Tight contact structures on fibered hyperbolic 3-manifolds}

\author{Ko Honda}
\address{University of Southern California, Los Angeles, CA 90089} 
\email{khonda@math.usc.edu}
\urladdr{http://math.usc.edu/\char126 khonda}

\author{William H. Kazez}
\address{University of Georgia, Athens, GA 30602}
\email{will@math.uga.edu}
\urladdr{http://www.math.uga.edu/\char126 will}

\author{Gordana Mati\'c}
\address{University of Georgia, Athens, GA 30602}
\email{gordana@math.uga.edu}
\urladdr{http://www.math.uga.edu/\char126 gordana}

\date{October 2, 2001.}

\thanks{KH supported by NSF grant DMS-0072853, AIM, and IHES;  GM supported by 
NSF grant DMS-0072853; WHK supported by NSF grant DMS-0073029.}

\begin{abstract}
We take a first step towards understanding the relationship between foliations 
and universally tight contact structures on hyperbolic 3-manifolds.  If  a 
surface bundle over a circle has pseudo-Anosov holonomy, we obtain a 
classification of  ``extremal'' tight contact structures.  Specifically, there 
is exactly one contact structure whose Euler class, when evaluated on the fiber, 
equals the Euler number of the fiber.  This rigidity theorem is a consequence of 
properties of the action of pseudo-Anosov maps on the complex of curves of the 
fiber and a remarkable flexibility property of convex surfaces in such a space. 
Indeed this flexibility may be seen in surface bundles over an interval where 
the analogous classification theorem is also established. 
\end{abstract}

\maketitle

Recent work of Colin \cite{Co3} and the authors \cite{HKM2}, and work in 
progress of Colin, Giroux, and the first author \cite{CGH} combine to suggest a 
very general, but rough, classification principle for universally tight contact 
structures.

\s\n
{\it Principle:} If $M$ is a closed, oriented, irreducible 3-manifold, then $M$ 
carries finitely many (isotopy classes of) universally tight contact structures 
if and only if $M$ is atoroidal. \s

The work \cite{CGH} on the finiteness of tight contact structures does not give 
an accurate bound on the number of universally tight contact structures carried 
by an atoroidal manifold.  (For example, it does not address the existence 
question.)  Indeed, it is not clear, {\it a priori}, just what are the 
properties of a 3-manifold which determine the number and nature of the tight 
contact structures it carries.  

As a first step towards understanding the classification of tight contact 
structures for atoroidal manifolds with infinite fundamental group, we consider 
hyperbolic 3-manifolds which fiber over the circle. These manifolds have enough 
structure to be accessible through a variety of techniques, and the 
relationships between the various structures they support are often extremely 
good predictors of similar relationships in more general 3-manifolds.  The 
uniqueness theorem (Theorem~\ref{thm:fibration}) suggests that contact topology 
may ultimately be a discrete version of foliation theory. More specifically, 
Eliashberg and Thurston \cite{ET} showed that any taut foliation may be 
perturbed to a universally tight contact structure.  There are many foliations 
of a fibered hyperbolic manifold that contain a fixed fiber as a leaf, yet 
perturbing any of them always produces the same universally tight contact 
structure by Theorem~\ref{thm:fibration}.

Contact structures on fibered hyperbolic 3-manifolds are studied by splitting 
the manifold along a fiber and later recreating the original manifold by gluing 
back using the pseudo-Anosov monodromy map.  This leads to an analysis of tight 
contact structures on the product $\Sigma\times [0,1]$, where $\Sigma$ is a 
closed, oriented surface of genus $g>1$.  There is a surprisingly small number 
of (isotopy classes of) tight contact structures on $\Sigma\times I$ with 
boundary conditions given in Theorem~\ref{classification}, and the theory for 
$g>1$ is quite different from the $g=1$ case treated in \cite{Gi00a, H1}.  Much 
of the classification on $\Sigma \times I$, at least in the cases we are led to 
study, can be encoded in the curve complex on $\Sigma$.  Therefore, 
Theorem~\ref{thm:fibration} exploits the interplay between the pseudo-Anosov 
monodromy and the curve complex of $\Sigma$.  The fact that there are only a few 
distinct isotopy classes of tight contact structures on $\Sigma \times I$ (of 
the kind we are interested in) is in large part due to a remarkable flexibility 
property (Proposition~\ref{ubiquity}) enjoyed by surfaces in $\Sigma \times I$ 
which are isotopic to $\Sigma\times\{t\}$.  To paraphrase the result, by 
isotoping such a surface, we are free to choose any pair of parallel, 
nonseparating curves as its dividing set.

\s\s

Let $\Sigma$ be a closed, oriented surface of genus $g>1$.  We study tight 
contact structures $\xi$ on the 3-manifold $M=\Sigma\times I=\Sigma\times[0,1]$ 
which satisfy the following condition:

\s\n
{\bf Extremal condition.}  $\langle 
e(\xi),\Sigma_t\rangle=\pm (2g-2)$, $t\in[0,1]$, where the left-hand side refers 
to the Euler class of $\xi$ evaluated on $\Sigma_t$.  

\s\n
Here we write $\Sigma_t=\Sigma\times\{t\}$. Tight contact structures on 
$\Sigma\times I$ satisfying this condition are said to be {\it extremal} for the 
following reason:  if $\xi$ is a tight contact structure on $\Sigma\times I$, 
then the Bennequin inequality \cite{Be,E92} states that:  
$$-(2g-2)\leq \langle e(\xi),\Sigma_t\rangle\leq 2g-2.$$ 

One of the main results of this paper is the following classification:

\begin{thm}   \label{classification}
Let $\Sigma$ be a closed oriented surface of genus $\geq 2$ and $M=\Sigma\times 
I$.  Fix dividing sets $\Gamma_{\Sigma_i}=2\gamma_i$ (2 parallel disjoint copies 
of $\gamma_i$), $i=0,1$, so that $\gamma_0$, $\gamma_1$ are nonseparating curves 
and $\chi((\Sigma_0)_+) - \chi((\Sigma_0)_-) = \chi((\Sigma_1)_+) - 
\chi((\Sigma_1)_-)$, where $(\Sigma_i)_+$, $(\Sigma_i)_-$ are the positive and 
negative regions of $\Sigma_i\setminus \Gamma_{\Sigma_i}$.   Then choose a 
characteristic foliation $\mathcal{F}$ on $\bdry M$ which is adapted to 
$\Gamma_{\Sigma_0}\sqcup \Gamma_{\Sigma_1}$.  We have the following:

\be
\item All the tight contact structures which satisfy the boundary condition 
$\mathcal{F}$ are universally tight. 

\item If $\gamma_0\not=\gamma_1$, then $\#\pi_0(Tight(M,\mathcal{F}))=4$.  
Provided $\gamma_0$ and $\gamma_1$ are not homologous, they are 
distinguished by the relative Euler classes $PD(\tilde{e}(\xi))=\pm\gamma_0\pm 
\gamma_1$. 

\item If $\gamma_0=\gamma_1$, then $\#\pi_0(Tight(M,\mathcal{F}))=5$. Three of 
them (one of them $I$-invariant) have relative Euler class 
$PD(\tilde{e}(\xi))=0$ and the other two have $PD(\tilde{e}(\xi))=\pm 
2\gamma_0=\pm 2\gamma_1$.

\ee

\end{thm}

Here, $\#\pi_0(Tight(M,\mathcal{F}))$ denotes the number of connected 
components of tight contact 2-plane fields adapted to $\mathcal{F}$, and the 
{\it relative Euler class} is an invariant of the tight contact structure which 
will be defined in Section~\ref{relative-euler}.    Theorem~\ref{classification} 
is a complete classification in the extremal case, provided 
$\#\Gamma_{\Sigma_i}=2$ and $\Gamma_{\Sigma_i}$ consists of two nonseparating 
curves.  The proof of Theorem~\ref{classification}, given in 
Section~\ref{section:classification} requires one involved calculation, found in 
Section~\ref{section:basecase}, followed by judicious use of general facts from 
curve complex theory, found in Section~\ref{section:curvecomplex}.

The extremal case is currently the only case we understand, but we also have the 
following theorem:

\begin{thm}   \label{thm:fibration}
Let $M$ be a closed, oriented, hyperbolic 3-manifold which fibers over $S^1$, 
where the fiber is a closed oriented surface $\Sigma$ of genus $g>1$ and the 
monodromy map is pseudo-Anosov.   Then there exists a unique tight contact 
structure in each of the two {\em extremal cases}, i.e., $\langle 
e(\xi),\Sigma\rangle = \pm(2g-2)$.  This contact structure is universally tight 
and weakly symplectically fillable.  Moreover, every $C^0$-small perturbation of 
the fibration into a contact structure is isotopic to the unique extremal tight 
contact structure.   
\end{thm}

Perturbing the fibration by $\Sigma$ into a contact structure is either done 
directly or by appealing to the perturbation result of Eliashberg and Thurston 
\cite{ET}, which also tells us that the contact structure is weakly 
symplectically semi-fillable.  It is interesting to note that, no matter how we 
perturb the fibration into a contact structure, the resulting tight contact 
structure is the same.  This contrasts with the case where the fiber is a torus 
\cite{Gi99a}.  

\s\n
In this paper we adopt the following conventions:

\be

\item The ambient manifold $M$ is an oriented, compact $3$-manifold.

\item $\xi$ = positive contact structure which is co-oriented by a global 1-form 
$\alpha$.

\item A convex surface $\Sigma$  is either closed or compact with Legendrian 
boundary.

\item $\Gamma_\Sigma$ = dividing multicurve of a convex surface $\Sigma$.

\item $\#\Gamma_\Sigma$ = number of connected components of $\Gamma_\Sigma$.

\item $\Sigma\setminus \Gamma_\Sigma=\Sigma_+\cup \Sigma_-$, where $\Sigma_+$  
(resp.\ $\Sigma_-$) is the region where the normal orientation of $\Sigma$ is 
the same as (resp.\ opposite to) the normal orientation for $\xi$. 

\item $|\beta\cap \gamma|$ = geometric intersection number of two curves $\beta$ 
and $\gamma$ on a surface.

\item $\#(\beta\cap\gamma)$ = cardinality of the intersection.

\item $\Sigma\setminus \gamma$ = metric closure of the complement of $\gamma$ in 
$\Sigma$.

\item $t(\beta,Fr_S)$ = twisting number of a Legendrian curve with respect to 
the framing induced from the surface $S$.

\ee

\section{Tools from convex surface theory}  

In this section we collect some results from the theory of convex surfaces which 
are nonstandard.  For standard results on convex surfaces, we refer the reader 
to \cite{Gi91,Gi00b,H1,HKM1,HKM2}.  

We first recall the Legendrian realization principle (LeRP), in a slightly 
stronger form.  An embedded graph $C$ on a convex surface $\Sigma$ is   
{\it nonisolating} if (1) $C$ is transverse to $\Gamma_\Sigma$, (2) the 
univalent vertices of $C$ lie on $\Gamma_\Sigma$, (3) all the other vertices 
do not lie on $\Gamma_\Sigma$, and (4) every component of $\Sigma\backslash 
(\Gamma_\Sigma \cup C)$ has a boundary component which intersects 
$\Gamma_\Sigma$.  

\begin{thm}[Legendrian realization] \label{lerp}
Let $C$ be a nonisolating graph on a convex surface $\Sigma$ and $v$ a contact 
vector field transverse to $\Sigma$. Then there exists an isotopy 
$\phi_s$, $s\in[0,1]$ so that: 
\be
\item $\phi_0=id$, $\phi_s|_{\Gamma_\Sigma}=id$.
\item $\phi_s(\Sigma)\pitchfork v$ (and hence $\phi_s(\Sigma)$ are all convex),
\item $\phi_1(\Gamma_\Sigma)=\Gamma_{\phi_1(\Sigma)}$,
\item $\phi_1(C)$ is Legendrian.
\ee
\end{thm}

\begin{prop}[The Right-to-Life Principle]  Let $S$ be a convex surface in a 
tight contact manifold $(M,\xi)$ and let $\delta$ be a Legendrian arc transverse 
to $\Gamma_S$ with endpoints on $\Gamma_S$, for which $|\Gamma_S\cap \delta|=3$.  
Suppose an ``abstract bypass move'' was applied to $\Gamma_S$ along $\delta$ and 
yielded a dividing set isotopic to $\Gamma_S$, i.e., this move was {\it 
trivial}.  (Here, by an ``abstract bypass move'' we simply mean a modification 
of the multicurve $\Gamma_S\subset S$ which would theoretically arise from 
attaching a bypass along $\delta$ --- there may or may not be an actual bypass 
half-disk along $\delta$ inside the ambient 3-manifold.)  Then there exists an 
actual bypass half-disk $B$ along $\delta$ (from the proper side) contained in 
the $I$-invariant neighborhood of $S$. \end{prop}

\n
The Right-to-Life Principle is a consequence of Eliashberg's classification of 
tight contact structures on the 3-ball \cite{E92}, and is proved in Lemma~1.8 
of \cite{H3}.  Here we recreate the proof.

\begin{proof}
Suppose $\delta$ intersects $\Gamma_S$ successively along $p_1,p_2,p_3$.  Since 
$\delta$ gives rise to a trivial ``abstract bypass move'', we may assume that 
the subarc from $p_1$ to $p_2$, together with a subarc of $\Gamma_S$, bounds a 
disk $D'$.  Let $D\subset S$ be a disk which contains $D'\cup \delta$.  We may 
assume $\bdry D$ is Legendrian with $tb(\bdry D)=-2$ --- to do this we apply 
LeRP (or Giroux's Flexibility Theorem \cite{Gi91,H1}) to $\bdry D$, while fixing 
$\delta$.     Then $\Gamma_D$ consists of two arcs $\gamma_1\ni p_1,p_2$ and 
$\gamma_2\ni p_3$.  

We will now find the actual bypass along $\delta\times\{0\}$ inside the 
$I$-invariant tight contact structure on $D\times [0,1]$.   (Here we are 
assuming that the bypass attached from the $D\times[0,1]$-direction is the 
trivial bypass.) Let $p_0$ be a point on $\gamma_2$ ($\not=p_3$) for which there 
exists a Legendrian arc $\delta'\subset D$ from $p_1$ to $p_0$ which does not 
intersect $\Gamma_D$ except at endpoints and does not intersect $\delta$ except 
at $p_1$.  Define $\delta_0 = \delta \cup \delta'$ and $\delta_1\subset D$ 
to be a Legendrian arc from $p_3$ to $p_0$ that lies on the same side of 
$\gamma_2$ as $\delta_0$ and has no other intersections with $\Gamma_D$.  
Now consider an arc $\varepsilon\subset D$ from $p_0$ to $p_3$ that lies on the 
opposite side of $\gamma_2$ as $\delta_0$ and $\delta_1$ and has no other 
intersections with $\Gamma_D$.  Let $A \subset D \times [0,1]$ be a convex 
annulus such that $\partial A = (\delta_0 \cup \varepsilon) \times \{ 0 \} \cup 
(\delta_1 \cup \varepsilon) \times \{1\}$ and such that $(\varepsilon \times 
[0,1]) \subset A$.  Since $\Gamma_{\varepsilon \times I}$ consists of an arc 
$\{q\}\times[0,1]$ ($q\in \varepsilon$), $\Gamma_A$ must contain one of two 
possible bypasses, one intersecting $\{p_0,p_1,p_2\}\times \{0\}$ and the other 
intersecting $\{p_1,p_2,p_3\}\times \{0\}$.  The former is a bypass which cannot 
exist since it gives rise to an overtwisted disk.  Therefore we obtain the 
second, as stated in the proposition.  \end{proof}

\begin{lemma}[Bypass Sliding Lemma]   \label{bypass-sliding}
Let $R$ be an embedded rectangle with consecutive sides $a, b, c, d$ 
in a convex surface $S$ such that $a$ is an arc of attachment of a 
bypass, $b$ and $d$ are subsets of $\Gamma_S$ and $c$ is a Legendrian 
arc which is efficient (rel endpoints) with respect to $\Gamma_S$.  Then there 
exists a bypass for which $c$ is its arc of attachment.
\end{lemma} 

\begin{proof}
The appendix to \cite{H2} explains how to explicitly move the endpoints of the 
Legendrian arc of attachment.  In this paper, we will deduce the Bypass Sliding 
Lemma from the Right-to-Life Principle.

Let $R'\supset R$ be an embedded disk in $S$ satisfying the following:

\s

\begin{itemize}
\item $\bdry R'$ is Legendrian, $\bdry R'\pitchfork \Gamma_S$, and $tb(\bdry 
R')=-3$. 

\item $\bdry R'=a'\cup b'\cup c'\cup d'$, where $a'$ and $c'$ are 
Legendrian arcs parallel to and close to $a$ and $c$, respectively.   The four 
arcs $a', b', c', d'$ are consecutive sides of a rectangle whose corners, lying 
on $\Gamma_S$,  have been smoothed. 

\item $\Gamma_{R'}$ consists of three parallel (= nested) dividing arcs.  

\end{itemize}

\s\n
Such an $R'$ exists by LeRP.  Suppose we attach the bypass along $a$ to obtain a 
convex surface $R''$ isotopic to $R'$ rel boundary.  Note that, away from $a$, 
$R''$ and $R'$ are identical.  Now, the ``abstract bypass move'' along $c$ 
applied to $R''$ still yields the same dividing set $\Gamma_{R''}$.  Therefore, 
the bypass must also exist along $c$.  (In fact, it exists in a small 
neighborhood of $R''$.) 
\end{proof}

\begin{lemma}\label{lemma:super-lerp}    Let $S$ be a closed convex surface of 
genus $g$ and $\Gamma_S$ its dividing set, consisting of one homotopically 
nontrivial separating curve.  Let $S\times I$ be its $I$-invariant neighborhood.  
Then there exists a bypass in $S\times I$ along $S\times\{1\}$ such that the 
dividing set $\Gamma_{S'}$ of the surface $S'$, obtained after bypass 
attachment, consists of three curves parallel to $\Gamma_S$.  Moreover, the arc 
of attachment $\alpha$ is contained in an annular neighborhood of $\Gamma_S$ and 
intersects $\Gamma_S$ in 3 points, and the two half disks bounded by $\alpha$ 
and $\Gamma_S$ have a common intersection along an arc contained in $\Gamma_S$.
\end{lemma}

For a more thorough discussion on dividing-curve increasing bypasses, see 
\cite{H3}.

\s\n
{\it Remarks.}  
\be

\item Typically, we increase $\#\Gamma$ by {\it folding} along a Legendrian  
divide (see \cite{H1} or  \cite{HKM2}).  This generates a pair of dividing 
curves parallel to the Legendrian divide.  However, this folding operation to 
increase $\#\Gamma$ does not occur immediately for free, in case the curve we 
want to realize as a Legendrian divide is an {\it isolating} curve for $\Gamma$ 
in the sense of LeRP.   This is the case in Lemma~\ref{lemma:super-lerp}.

\item Lemma~\ref{lemma:super-lerp} can be viewed as a strengthened form of LeRP 
and of Giroux's Flexibility Theorem.

\item Lemma~\ref{lemma:super-lerp} (or the proof therein) shows that folds along 
homotopically nontrivial closed curves always exist.

\ee

\begin{proof} In order to use LeRP, we first cut along some annulus 
$\gamma^*\times I$, where $\gamma^*\subset S$ is a closed nonseparating curve 
which does not intersect $\Gamma_S$.   Let $N$ be $(S\times I)\setminus 
(\gamma^*\times I)$, with edges rounded near $\gamma^*\times I$, so that $\bdry 
N$ is convex.  We will write $G=\bdry N$. Note that we may think of  
$\Gamma_S$ as a sub-multicurve of  $\Gamma_G$.  A parallel copy $(\Gamma_S)^*$ 
of $\Gamma_S$ on $G$, pushed off closer towards $\gamma^*\times I$ and disjoint 
from $\Gamma_S$, is then {\it nonisolating}.  We can then use LeRP to realize 
$(\Gamma_S)^*$ as a Legendrian divide and use it to perform a fold. (For more 
details on folding and bypasses, see \cite{HKM2}.   In a standard $I$-invariant 
neighborhood of $G$, we therefore find a parallel copy $G'$, where $\Gamma_{G'}$ 
is $\Gamma_G$, with 2 parallel curves (parallel to $\Gamma_S$) added.  To get 
from $G'$ to $G$, a bypass $B$ is attached along $G'$ which straddles the three 
distinct dividing curves parallel to $\Gamma_S$.  This has the effect of 
reducing $\#\Gamma$ by $2$.  Every bypass operation always has its inverse 
operation, and the inverse of $B$ is $B^{-1}$, which is attached along $G$ and 
increases $\#\Gamma$ by $2$.  We can now view $B^{-1}$ as also attached onto $S$ 
and intersecting $\Gamma_S$ three times.  The Legendrian arc of attachment for 
$G$ may not exist on $S$, but one can always find an appropriate arc of 
attachment using Giroux's Flexibility Theorem. (Remark:  The arc of attachment 
does not really matter here --- the only thing that really matters is the 
twisting number of the bypass Legendrian arc.)  Hence the bypass survives 
passage from being a bypass for $\Gamma_G$ to being a bypass for $\Gamma_S$.   
\end{proof}

\section{Tools from Curve Complex Theory} \label{section:curvecomplex}

In this section, we recall several facts from the theory of curve complexes 
which will become useful later.  Let $\Sigma$ be a closed oriented surface of 
genus $g\geq 2$  --- we stress here that all the lemmas and propositions of this 
section require $g\geq 2$.   Let $\mathcal{C}(\Sigma)$ be the curve 
complex for $\Sigma$. The vertices of $\mathcal{C}(\Sigma)$ consist of isotopy 
classes of homotopically nontrivial closed curves on $\Sigma$ and there is a 
single edge between each pair $\gamma_0$, $\gamma_1$ of (isotopy classes of) 
closed curves with $|\gamma_0\cap \gamma_1|=0$.  Note that we do not attach 
simplices of higher dimension, since they are not needed. The facts we need 
from the theory of curve complexes are variations and strengthenings of the 
basic Lemma~\ref{basiclemma} below.  We refer the reader to \cite{Iv} for an 
exposition on curve complexes, an extensive reference, and a proof of 
Lemma~\ref{basiclemma}.  

\begin{lemma} \label{basiclemma}
$\mathcal{C}(\Sigma)$ is connected, i.e., given two closed curves $\alpha$, 
$\alpha'$ in $\mathcal{C}(\Sigma)$,  there exists a 
sequence $\alpha_0=\alpha,\alpha_1,\dots,\alpha_k=\alpha'$ in 
$\mathcal{C}(\Sigma)$ where $|\alpha_{i-1}\cap\alpha_i|=0$, $i=1,\dots,k$.
\end{lemma}

Alternatively, we have the following:

\begin{lemma}\label{curvecomplex}
Given two nonseparating closed curves $\alpha$, $\alpha'$ on $\Sigma$, there
exists a sequence $\alpha_0=\alpha,\alpha_1,\dots,\alpha_k=\alpha'$ of closed
curves where $|\alpha_{i-1}\cap\alpha_i|=1$, $i=1,\dots,k$.
\end{lemma}

\begin{prop}\label{fact1}
Given two nonseparating closed curves $\alpha$ and $\alpha'$ on 
$\Sigma$, there exists a sequence $\alpha_0=\alpha, \alpha_1, \dots,
\alpha_k=\alpha'$ of closed curves where:

\be
\item $\alpha_i$ are nonseparating, $i=0,\dots,k$.
\item $|\alpha_{i-1}\cap \alpha_{i}|=0$, $i=1,\dots, k$.
\item $\alpha_{i-1}$ and $\alpha_{i}$ are not homologically equivalent, 
$i=1,\dots, k$. 
\ee
\end{prop}

\begin{proof}
By Lemma~\ref{curvecomplex}, there is a sequence 
$\alpha_0=\alpha,\alpha_1,\dots,\alpha_k=\alpha'$ of closed  curves where 
$|\alpha_{i-1}\cap\alpha_i|=1$, $i=1,\dots,k$.  A neighborhood $U_{i-1}$  of 
$\alpha_{i-1} \cup \alpha_i$ is a punctured torus.  Since the genus of $\Sigma$ 
is at least 2, $\Sigma - U_{i-1}$ contains a nonseparating curve $\beta_{i-1}$.  
Therefore, $\alpha_0, \beta_0, \alpha_1, \beta_1, \dots, \alpha_k$ is the 
desired sequence.
\end{proof}

\begin{prop} \label{fact2}
Given three nonseparating closed curves $\alpha$, $\beta$, $\beta'$ 
on $\Sigma$, where $|\alpha\cap \beta|=0$, $|\alpha\cap\beta'|=0$, $\alpha$ and 
$\beta$ are nonhomologous, and $\alpha$ and $\beta'$ are nonhomologous, there 
exists a sequence $\beta_0=\beta, \beta_1, \dots, \beta_k=\beta'$ of closed 
curves where:

\be
\item $\beta_i$ are nonseparating, $i=0,\dots,k$.
\item $|\alpha\cap\beta_i|=0$, $i=0,\dots, k$.
\item $|\beta_{i-1}\cap \beta_i|=1$, $i=1,\dots, k$.
\ee

\end{prop}

\begin{proof}
By assumption $\Sigma \setminus \alpha$ is connected.  If $\gamma$ is a curve in 
$\Sigma$ such that $|\gamma \cap \alpha|=0$, then $\gamma$ is nonseparating and 
not homologous to $\alpha$ if and only if $\gamma$ is nonseparating as a subset 
of $\Sigma \setminus \alpha$.  Let $\Sigma'$ be $\Sigma \setminus \alpha$ 
together with two disks capping off the boundary components corresponding to 
$\alpha$. Then $\beta$ and $\beta'$ are nonseparating closed curves of $\Sigma'$ 
and applying Lemma~\ref{curvecomplex} to $\beta, \beta'$ and $\Sigma'$ produces 
the desired sequence of curves.
\end{proof}

\begin{prop} \label{fact0}
Given two nonseparating closed curves $\alpha$ and $\alpha'$ on $\Sigma$, there 
exists a sequence $\alpha_0=\alpha, \alpha_1, \dots, \alpha_k=\alpha'$ of closed 
curves where:

\be
\item $\alpha_i$ are nonseparating, $i=0,\dots k$.
\item $|\alpha_{i-1}\cap \alpha_{i}|=1$, $i=1,\dots,k$.
\item $|\alpha_{i-1}\cap \alpha_{i+1}|=0$, $i=1,\dots,k-1$.
\ee

\end{prop}

\begin{proof}
Let $\alpha_0=\alpha, \alpha_1, \dots, \alpha_k=\alpha'$ be a sequence for which
$|\alpha_{i-1} \cap \alpha_i| =1$, $i=1,\dots,k$, as guaranteed by 
Lemma~\ref{curvecomplex}. Suppose that $j \ge 1$ is the smallest integer for 
which $|\alpha_{j-1} \cap \alpha_{j+1}| > 0$. We show how to insert curves into 
the sequence between $\alpha_j$ and $\alpha_{j+1}$, until Condition~3 is 
satisfied, at least up to $\alpha_{j+1}$.  Since curves are only inserted after 
$\alpha_{j}$, there is no danger that Condition~3 will be disturbed for earlier 
terms in the sequence.

For notational convenience, we shall take $j=1$. We will always assume
that $\alpha_0 \cap \alpha_1$ is not an element of $\alpha_2$. Suppose that
$|\alpha_0 \cap \alpha_2| \ge 1$.  The following cases are used to inductively
decrease the number of intersections between $\alpha_0$ and $\alpha_2$.

\s\n
{\bf Case~1:} Not all intersections between $\alpha_0$ and $\alpha_2$ have 
the same sign.

\s\n
Then there exist two points $p, q\in \alpha_0\cap\alpha_2$ with opposite signs 
of intersection, and an arc $[p,q]_{\alpha_0}\subset \alpha_0$ which 
connects $p$ and $q$ such that $[p,q]_{\alpha_0}$ does not intersect $\alpha_1$ 
and contains no point of $\alpha_2$ in its interior.  Cut-and-paste surgery of 
$\alpha_2$ along $[p,q]_{\alpha_0}$ produces two curves; let $\beta$ the 
component which intersects $\alpha_1$ once.  The sequence $\alpha_0, \alpha_1, 
\beta, \alpha_1, \alpha_2$ satisfies Condition~2 (and hence Condition~1), and is 
closer to satisfying Condition~3 than the original sequence, in the sense that 
the total number of intersection points of type $\alpha_{i-1}\cap 
\alpha_{i+1}$ has been reduced.  Here it is understood that the second 
$\alpha_1$ in the sequence is slightly isotoped off the first $\alpha_1$ so that 
the two curves are disjoint as required in Condition~3.

\s\n
{\bf Case~2:}  All intersections of $\alpha_0$ and $\alpha_2$ have the same sign 
and $|\alpha_0 \cap \alpha_2| \ge 2$.

\begin{figure}[ht]
	{\epsfysize=106pt\centerline{\epsfbox{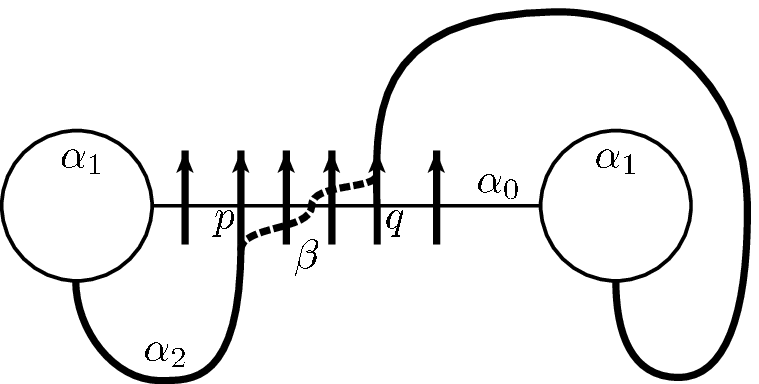}}}
	\caption{Case 2}
	\label{Case2}
\end{figure}

\s\n
Let $[p,q]_{\alpha_2}$, $p,q\in \alpha_0\cap \alpha_2$, be the smallest interval 
in $\alpha_2$ that contains $\alpha_1 \cap \alpha_2$.  
Let $[p,q]_{\alpha_0}$ be the interval in $\alpha_0$ (with endpoints $p,q$) 
which does not intersect $\alpha_1$.  Let $\beta = [p,q]_{\alpha_2} \cup 
[p,q]_{\alpha_0}$. After a slight isotopy of $\beta$, the sequence $\alpha_0, 
\alpha_1, \beta, \alpha_1,\alpha_2$ satisfies Condition~2. Condition~3 is not 
satisfied, but at least the first and third terms of this interim sequence 
intersect only once and $|\beta \cap \alpha_2| < |\alpha_0 \cap
\alpha_2|$.

\s\n
{\bf Case~3:}  $|\alpha_0 \cap \alpha_2| = 1$.

\s\n
Let $U$ be a regular neighborhood of $\alpha_0 \cup \alpha_1 \cup \alpha_2$ in
$\Sigma$.  Let $V = U \setminus (\alpha_0 \cup \alpha_1)$.  $V$ is a planar 
surface with 4 boundary components, one corresponding to $\alpha_0 \cup 
\alpha_1$, shown as a rectangle in Figure~\ref{Case3}, and three others denoted 
$A, B$ and $C$.  The construction of the desired sequence of curves depends on 
the relation of $A, B$ and $C$ to the rest of $\Sigma$.  In each of the 
following cases, $S$ will denote a connected subsurface of $\Sigma$ which 
satisfies $\partial S \subset A \cup B \cup C$ and whose interior $int(S)$ is 
disjoint from $V$.

\begin{figure}[ht]
	{\epsfysize=4.5in\centerline{\epsfbox{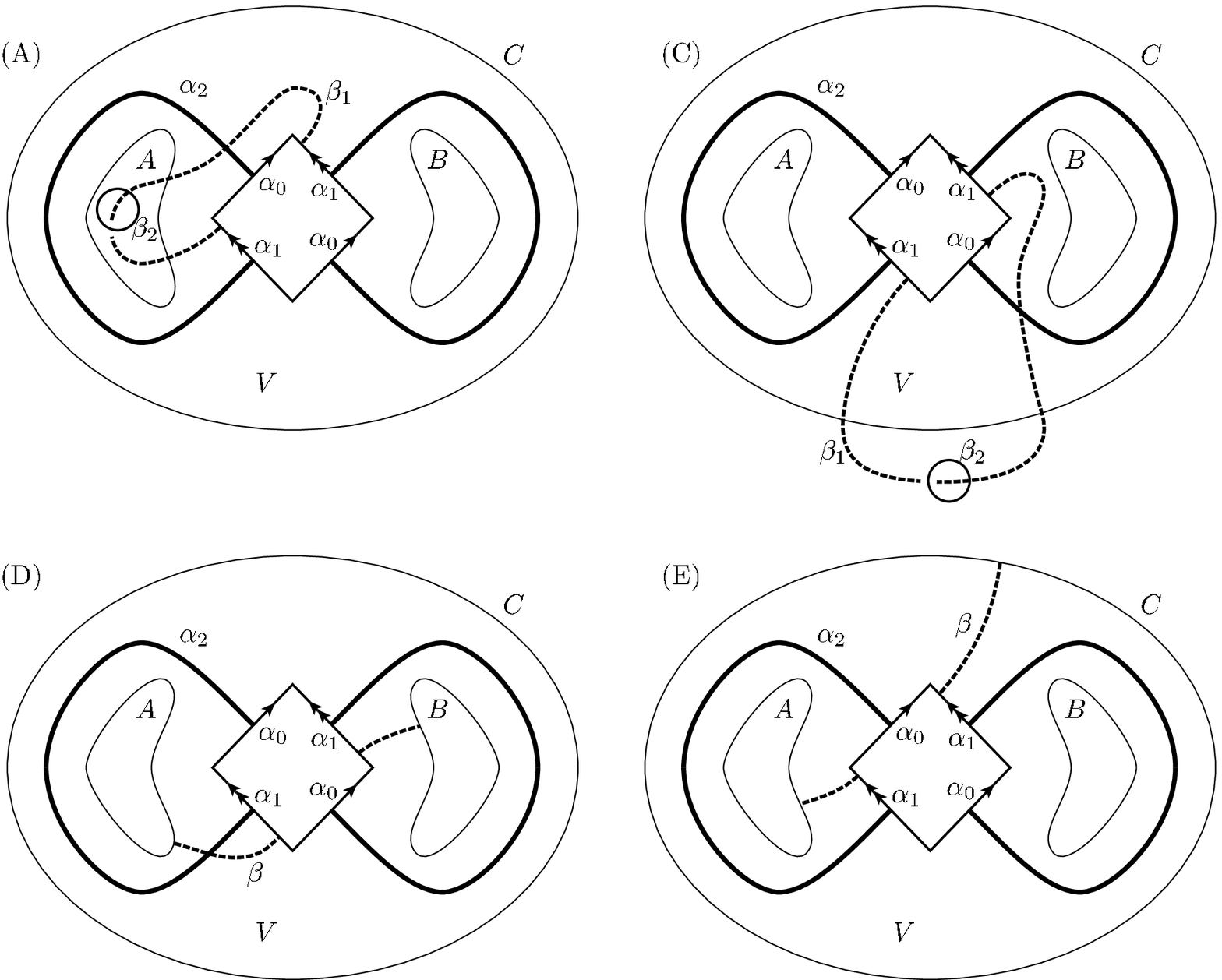}}}
	\caption{Case 3}
	\label{Case3}
\end{figure}

\s\n
(A) $A\subset S$ and $\text{genus}(S) \ge 1$.  Choose $\beta_1$ as indicated in
Figure~\ref{Case3}(A), that is, $\beta_1$ starts on $\alpha_1$, then 
passes across
$A$ into $S$ and back to $A$ without separating $S$, and then crosses 
$\alpha_2$
before returning to the point it started from, but from the opposite side of
$\alpha_1$.  Let $\beta_2 \subset S$ be a curve such that $|\beta_2 \cap
\beta_1| = 1$.  The sequence $\alpha_0, \alpha_1, \beta_1, \beta_2, \beta_1,
\alpha_2$ satisfies Conditions~1--3.

\s\n
(B) $B\subset S$ and $\text{genus}(S) \ge 1$.  This is similar to (A).

\s\n
(C) $C\subset S$ and $\text{genus}(S) \ge 1$.  Choose $\beta_1$ and 
$\beta_2$ as
shown in Figure~\ref{Case3}(C) and the sequence $\alpha_0, \alpha_1, \beta_1,
\beta_2, \beta_1, \alpha_2$ will satisfy Conditions~1--3.

\s\n
(D) $A \cup B \subset S$.  Choose an arc in $V$ that runs from $A$, across
$\alpha_2$ and $\alpha_1$ and then to $B$.  Complete this arc to a 
curve $\beta$
by attaching an arc in $S$.  The sequence $\alpha_0, \alpha_1, \beta, A, \beta,
\alpha_2$ satisfies Conditions~1--3.

\s\n
(E) $A \cup C \subset S$.  Choose an arc in $V$ that runs from $C$, across 
$\alpha_1$ and then to $A$.  Complete this arc to a curve $\beta$ by attaching 
an arc in $S$.  The sequence $\alpha_0, \alpha_1, \beta, \alpha_1, \alpha_2$ 
satisfies Conditions~1--3.

\s\n
(F) $B \cup C \subset S$.  This is similar to (E).

\s\n
The only possibility not covered in (A)--(F) is if each of $A, B$ and $C$ bound
disks, but this would imply $\text{genus}(\Sigma) = 1$, contrary to assumption.
\end{proof}

\section{Definition of the relative Euler class}     \label{relative-euler}

Let $M=\Sigma\times I=\Sigma\times[0,1]$ and $\Sigma_t=\Sigma\times\{t\}$.  If
$\Sigma_t$ is convex, we denote $\Gamma_t=\Gamma_{\Sigma_t}$.  Let $\xi$ be a 
tight contact structure on $M$ which satisfies $\langle e(\xi),\Sigma_t 
\rangle=\pm (2g-2)$.   Since the situation is symmetric under sign change, we 
will additionally assume that $\langle e(\xi),\Sigma_t \rangle=-(2g-2)$.

First recall the following identity (cf. Kanda \cite{K98} or Eliashberg 
\cite{E92}): 
\begin{equation}      
\label{Eqn1} \langle e(\xi),S \rangle = \chi(S_+)-\chi(S_-).
\end{equation}
Here $S$ is a closed convex surface.  

According to Giroux's criterion \cite{Gi00b}, a closed convex surface 
$\Sigma\not=S^2$ has a tight neighborhood if and only if no connected component 
of $\Sigma_\pm$ is a disk.  This implies that $\chi(\Sigma_\pm)\leq 0$.  Hence, 
the extremal condition implies that $\chi(\Sigma_+)=-(2g-2)$ and  
$\chi(\Sigma_-)=0$.  In other words, $\Sigma_-$ is a nonempty union of annuli.  
Here, recall that both $\Sigma_+$ and $\Sigma_-$ of a convex surface $\Sigma$ 
are nonempty, since there must be both sources ($\Sigma_+$) and sinks 
($\Sigma_-$) for the characteristic foliation.  

This section is devoted to defining the {\it relative Euler class} 
$\tilde{e}(\xi)$ of $\xi$, not to be confused with the Euler class $e(\xi)$ used 
previously. It is sufficient to define the relative Euler class on annuli, as 
follows. Let $\gamma\subset \Sigma$ be a closed curve.  Suppose first that 
$\gamma\times\{0,1\}$ is {\it efficient} with respect to $\bdry M$, i.e., 
intersects $\Gamma_{\bdry M}$ minimally in its isotopy class in $\bdry M$.  If 
$\gamma\times\{0,1\}$ is {\it nonisolating} in $\bdry M$ (which is the case for 
example when $\gamma$ is nonseparating in $\Sigma$), then we may use the 
Legendrian realization principle (LeRP) \cite{H1} to make $\gamma\times\{0,1\}$ 
Legendrian.  If $S=\gamma\times I\subset M$ is a convex surface, then we define: 
\begin{equation}  \label{Eqn2}
\langle \tilde{e}(\xi),S\rangle = \chi(S_+) - \chi(S_-). 
\end{equation}
In general, choose $\gamma\times\{0,1\}$ (i.e., introduce extra intersections 
with $\Gamma_{\bdry M}$) so that the following condition 
holds:

\s\n
{\it $\gamma\times\{0,1\}$ is nonisolating and every 
connected component of $(\gamma\times\{i\})\cap (\Sigma_i)_-$, $i=0,1$, is a 
nonseparating arc.}

\s\n 
If this condition is satisfied, we will say that $\gamma\times\{0,1\}$ has been 
{\it primped}.  We remark here that (i) $(\Sigma_i)_-$ is a union of annuli, and 
(ii) components of $(\gamma\times \{i\})\cap (\Sigma_i)_+$ may be separating.  
We take $S$ convex with primped $\bdry S=\gamma\times\{0,1\}$ and define 
$\langle \tilde{e}(\xi),S\rangle$ as in Equation~\ref{Eqn2}.

What we would like to prove is that $\tilde{e}(\xi)$ is indeed a homology 
invariant, i.e.,  lives in $H^2(M,\bdry M;\Z)$.  This follows from the following 
three lemmas.

\begin{lemma}  \label{lemma1}
Let $S$ and $S'$ be convex surfaces with $\bdry S = \bdry S'$ primped 
Legendrian curves on $\bdry M$. If $S$ and $S'$ are isotopic rel boundary, then 
$\langle \tilde{e}(\xi),S\rangle = \langle \tilde{e}(\xi),S'\rangle$. 
\end{lemma}

\begin{proof}
Since $\bdry S=\bdry S'$, we consider $S\cup (-S')$.   After rounding along the 
common edge and perturbing slightly if necessary (without changing the 
isotopy class of the dividing set), we may take $S\cup (-S')$ to be a 
closed immersed surface.  Although Equation~\ref{Eqn1} holds for closed 
convex surfaces (convex surfaces are embedded by definition), we may clearly 
extend Kanda's argument in \cite{K98} to immersed surfaces which have meaningful 
positive and negative regions.  Therefore, we have:
\begin{eqnarray*}
\langle \tilde{e}(\xi), S\rangle -\langle\tilde{e}(\xi), S' \rangle & = &
[\chi(S_+)-\chi(S_-)]-[\chi(S_+')-\chi(S_-')] \\
&=& \chi(S_+\cup S'_-)-\chi(S_-\cup S'_+)   \\
& = & \langle e(\xi),S\cup(-S')\rangle = 0.
\end{eqnarray*}
This proves that $\langle \tilde{e}(\xi),S\rangle$ is independent of the choice 
of $S$, provided $\bdry S$ is fixed. \end{proof}

\begin{figure}[ht]
 {\epsfysize=120pt\centerline{\epsfbox{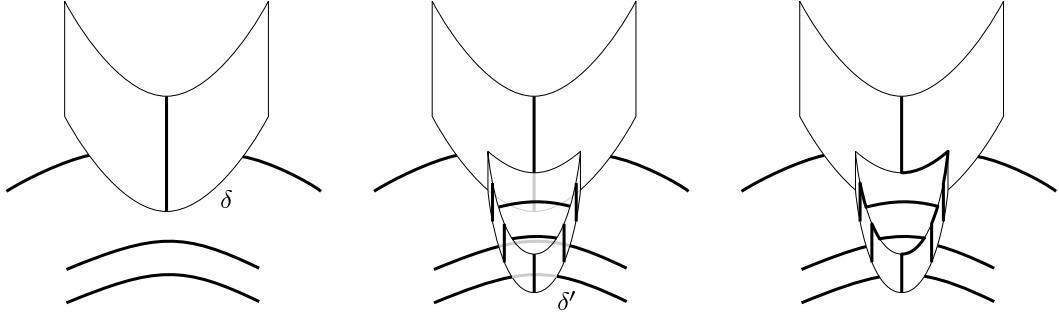}}}
 \caption{Isotopy of $\delta$ to $\delta'$}
 \label{primped.eps}
\end{figure}

\begin{lemma}  \label{extra-intersections}
Let $S=\delta\times I$ and $S'=\delta'\times I$ be isotopic convex surfaces with 
primped Legendrian boundary in $M$.  Then $\langle 
\tilde{e}(\xi),S\rangle=\langle \tilde{e}(\xi),S'\rangle$.
\end{lemma}

\begin{proof}
The key point of using primped curves is that dividing curves of 
$\Sigma\times\{0,1\}$ (in the extremal case) come in pairs, and we can isotop 
one primped curve to another primped curve in the same isotopy class by pushing
$\delta\times\{0,1\}$ across  two curves in a pair simultaneously, i.e. by 
changing the annulus by a sequence of operations of the type described in  
Figure~\ref{primped.eps} (or its inverse). We see that such an isotopy changes 
the dividing set on the annulus by adding (or deleting) a nested pair of 
boundary compressible arcs. This amounts to adding  (or removing) one disk 
region each for $S_+$ and $S_-$ and does not change $\langle 
\tilde{e}(\xi),S\rangle$.  A finite sequence of such steps will bring us to the 
situation where we can apply the previous lemma.
\end{proof}

Next, we prove additivity.  

\begin{lemma}
Let $\delta_i$, $i=1, 2$, be two oriented nonseparating curves on $\Sigma$ with 
nontrivial geometric intersection. Let $S_i=\delta_i\times I$ and form 
$S=S_1+S_2=\delta\times I$ by the natural cut-and-paste operation corresponding 
to homology addition.  Then $$\langle \tilde{e}(\xi),S\rangle   = \langle 
\tilde{e}(\xi),S_1\rangle + \langle \tilde{e}(\xi),S_2\rangle.$$ 
\end{lemma}

\begin{proof}
Consider the graph $C=\delta_1\cup \delta_2$, where  
$\delta_1\pitchfork \delta_2$ and $\delta_1\cap \delta_2$ does not intersect 
$\Gamma_\Sigma \cup \Sigma_-$.  Assume enough extra intersections of 
$\Gamma_\Sigma\cap \delta_i$ have been introduced to $\delta_i$, $i=1,2$, so 
that $\delta_i$ is primped and $C$ satisfies the {\it nonisolating} condition.   
We may take the common intersections $\delta_1\cap \delta_2$ to be elliptic 
tangencies, using a slightly stronger version of LeRP which is easily derived 
from Giroux's Flexibility Theorem \cite{H1}, and $(\delta_1\cap \delta_2)\times 
I$ to be transverse curves, by perturbing if necessary.  Let 
$\delta'=\delta_1+\delta_2$ be the multicurve obtained by smoothing the 
intersection $\delta_1\cap \delta_2$ in the standard manner.  Then the surface 
$S'$ obtained by performing a cut-and-paste along the transverse curve and 
smoothing the corners satisfies the following equality: $$\langle 
\tilde{e}(\xi),S'\rangle   = \langle \tilde{e}(\xi),S_1\rangle + \langle 
\tilde{e}(\xi),S_2\rangle.$$ The equality follows from relating 
$\chi(S'_+)-\chi(S'_-)$ to the more standard way of computing $\langle 
\tilde{e}(\xi),S'\rangle$ using signs and types of isolated singularities as in 
Kanda's argument in \cite{K98}.     Now, $\bdry S'$ is primped, since we took 
the intersections $\delta_1\cap \delta_2$ to be away from $\Sigma_-$. Therefore, 
Lemma~\ref{extra-intersections} implies that $\langle 
\tilde{e}(\xi),S'\rangle=\langle \tilde{e}(\xi),S\rangle$. \end{proof}

We will usually express $\tilde{e}(\xi)\in H^2(M,\bdry M;\Z)$ in terms of its 
Poincar\'e dual $PD(\tilde{e}(\xi))\in H_1(M;\Z)$.

\section{Computation of the base case}   \label{section:basecase}

Recall that, by our choice of $\langle e(\xi), \Sigma_t\rangle$, $(\Sigma_i)_-$, 
$i=0,1$, is a union of annuli. If we assume that there is only one pair of 
dividing curves on each  $\Sigma_i$,  then $(\Sigma_i)_-$ is an annulus and 
$(\Sigma_i)_+$ is its complement in $\Sigma_i$ (i.e., ``most'' of the surface).

We will now consider the following special case, which turns out to be the most 
fundamental:

\s\n
{\bf Base Case.}  $\Gamma_i=2\gamma_i$, $i=0,1$, and $|\gamma_0\cap 
\gamma_1|=1$.

\s
Here, $k\gamma$ is shorthand for $k$ parallel (mutually nonintersecting) 
copies of a closed curve $\gamma$.  

\begin{thm}      \label{basecase}
Let $\Sigma$ be a closed surface of genus $g\geq 2$, $M=\Sigma\times I$ the 
ambient manifold, and $\mathcal{F}$ a characteristic foliation on $\bdry M$ 
adapted to $\Gamma_0\sqcup \Gamma_1$ of the Base Case. 
Let $Tight(M,\mathcal{F})$ be the space of tight contact 2-plane fields which 
induce $\mathcal{F}$ along $\bdry M$. Then we have the following:

\be
\item $\#\pi_0(Tight(M,\mathcal{F}))=4$ .
\item The isotopy classes of tight contact structures are distinguished by their 
relative Euler class, which are: $$PD(\tilde{e}(\xi))=\pm \gamma_0\pm 
\gamma_1\in H_1(M;\Z).$$ \item All the tight contact structures are universally 
tight. \ee 

\end{thm}

This computation is a little involved, and occupies Sections~\ref{section21} 
through \ref{section24}.  In what follows, when we prescribe a boundary 
condition for a 3-manifold $M$, we will simply give $\Gamma_{\bdry M}$, 
although, strictly speaking, we need to also assign the characteristic foliation 
$\mathcal{F}$ on $\bdry M$ {\it adapted to} $\Gamma_{\bdry M}$.  We will assume 
that some convenient characteristic foliation is prescribed, since the actual 
number of tight contact structures is independent of the actual characteristic 
foliation adapted to $\Gamma_{\bdry M}$ (see \cite{H1}).    The following is a 
preliminary lemma.

\begin{lemma}  \label{uniqueness-bdry-case}
There exists a unique tight contact structure on $H=S\times I$, where 
$S$ is a compact oriented surface with nonempty boundary, 
$S_0$, $S_1$, and $(\bdry S)\times I$ are convex,  $\Gamma_{S_0}=\Gamma_{S_1}$ 
are $\bdry$-parallel, and $\Gamma_{\bdry S\times I}$ are vertical.  (Here {\it 
vertical} means the dividing curves are arcs $\{pt\}\times I$.)  This tight 
contact structure is universally tight, and is obtained by perturbing the 
foliation of $H$ by leaves $S\times \{t\}$, $t\in[0,1]$, into a contact 
structure. 
\end{lemma}

\begin{proof}
Observe that $H=S\times I$ is a handlebody and that $\Gamma_{\bdry H}$, after 
edge-rounding \cite{H1}, is isotopic to $\bdry S \times\{{1\over 2}\}$.   See 
Figure~\ref{fig1}.

\begin{figure}[ht] 	
	{\epsfysize=2.5in\centerline{\epsfbox{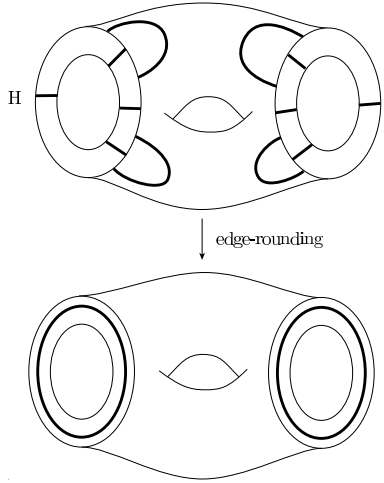}}}
	\caption{Dividing multicurve after edge-rounding.} 	
	\label{fig1}
\end{figure}

Consider a meridional disk $D=\gamma\times I$, where $\gamma$ is a  
non-boundary-parallel, properly embedded arc on $S$.  Then, after rounding, 
$\bdry D$ intersects $\Gamma_{\bdry H}$ along exactly two points.  Make $\bdry 
D$ Legendrian and $D$ convex.  Then the Thurston-Bennequin invariant $tb(\bdry 
D)$ equals $-1$, and there is a unique dividing curve configuration on $D$ 
consistent with this boundary condition.  Moreover, there is a system of such 
meridional disks $D_1,\dots,D_k$ which decompose $H$ into 3-balls $B^3$, each 
with $\Gamma_{\bdry B^3}=S^1$.  By Eliashberg's uniqueness theorem for tight 
contact structures on the 3-ball (c.f. \cite{E92}), there is a unique tight 
contact structure up to isotopy rel boundary on each of the $B^3$'s.  This 
implies that there exists at most one tight contact structure on $H$ with given 
$\Gamma_{\bdry H}$.  Now, to prove that there indeed exists a (universally) 
tight contact structure on $H$ with given $\Gamma_{\bdry H}$, we glue back using 
Theorem~\ref{colin-gluing} below. Note that $\Gamma_{D_i}$ is $\bdry$-parallel 
on each meridional disk $D_i$.  We leave the statement of the perturbation to 
the reader. \end{proof}

\begin{thm}[Colin~\cite{Co99a}] \label{colin-gluing}
Let $(M,\xi)$ be an oriented, compact, connected, irreducible, contact 
3-manifold and $S\subset M$ an incompressible convex surface with nonempty 
Legendrian boundary and $\bdry$-parallel dividing set $\Gamma_S$.  If  
$(M\setminus S, \xi|_{M\setminus S})$ is universally tight, then $(M,\xi)$ is 
universally tight. 
\end{thm}

\begin{figure}[ht]	
	{\epsfysize=2.5in\centerline{\epsfbox{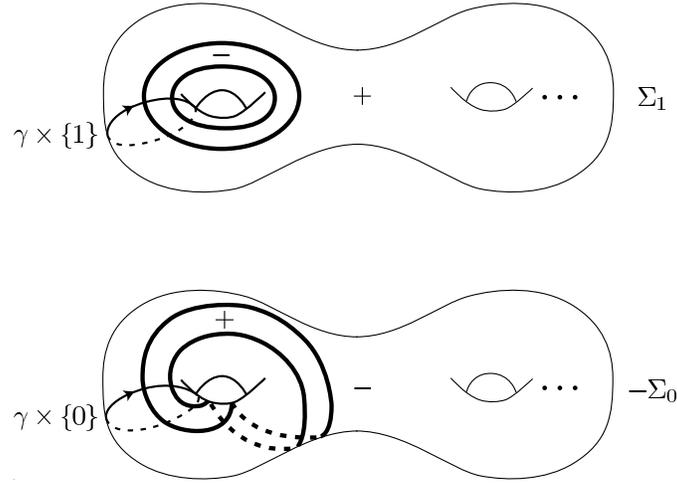}}}
	\caption{Suitable choice of $\Gamma_0$ and $\Gamma_1$.} 	
	\label{fig2}
\end{figure}

We now begin the analysis of the Base Case.  Let us position the 
dividing curves $2\gamma_0$ and $2\gamma_1$ as in Figure~\ref{fig2}, that is, we    
suppose $\#(\gamma_0\cap\gamma_1)=1$ and we have oriented $\gamma_0$ and 
$\gamma_1$ so that the intersection pairing $\langle \gamma_1,\gamma_0\rangle$ 
on $\Sigma$ equals $+1$. Now consider an oriented closed curve 
$\gamma$ which satisfies $\#(\gamma\cap\gamma_i)=1$ and 
$\langle\gamma,\gamma_i\rangle=1$, $i=0,1$.  For our convenience, we will assume 
that $\gamma_0$ and $\gamma_1$ are identical except in a thin annulus 
$\mathcal{A}$ parallel to but disjoint from $\gamma$, obtained as follows.  
First push $\gamma$ off of itself in the direction opposite the direction given 
by the orientation of $\gamma_1$ to obtain $\gamma'$.    Then let $\mathcal{A}$ 
be a small annular neighborhood of $\gamma'$.  $\gamma_0$ is then obtained from 
$\gamma_1$ via a Dehn twist in $\mathcal{A}$.  We will write 
$\Gamma_{\Sigma_1}\cap \gamma=\{p_1,p_2\}\times \{1\}$, $\Gamma_{\Sigma_0}\cap 
\gamma=\{p_1,p_2\}\times \{0\}$.  Our first cut of $M=\Sigma\times I$ will be 
along the convex surface $A=\gamma\times I$, where $A=\gamma\times I$ is given 
the orientation induced from $\gamma$ and $I$.

\begin{figure}[ht] 		
{\epsfysize=3in\centerline{\epsfbox{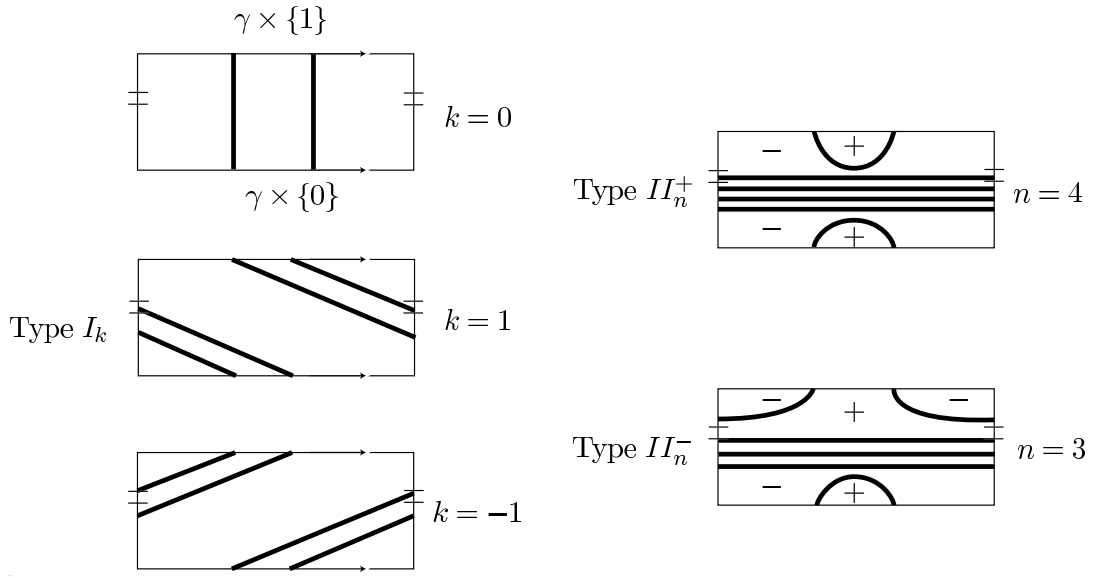}}} 	
\caption{Dividing curves of type ${\mbox{\it I}}_k$ and  $\ii_n^\pm$.} 		
\label{fig3} 
\end{figure}

There are two general classes of dividing curves $\Gamma_A$ which we denote by 
${\mbox{\it I}}_k$ and $\ii_n^\pm$.  The dividing set ${\mbox{\it I}}_k$ consists of 2 parallel 
nonseparating dividing curves (i.e., dividing curves which go across from 
$\gamma\times \{0\}$ to $\gamma\times \{1\}$).  Here, $k\in \Z$ denotes the {\it 
holonomy}, or the amount of spiraling, defined as follows.  First, {\it zero 
holonomy} $k=0$ means the dividing curves are {\it vertical} in the sense that 
they are isotopic rel boundary to $\{q_1,q_2\}\times [0,1]$, where 
$\{q_1,q_2\}\times\{0,1\}$ are the endpoints of $\Gamma_A$.  The holonomy is $k$ 
if $\Gamma_A$ is obtained from $\{q_1,q_2\}\times [0,1]$ by doing $k$ negative 
Dehn twists along the core curve of $A$. The dividing set $\ii_n^+$ (resp.\ 
$\ii_n^-$), $n\in \Z^{\geq 0}$, consists of two $\bdry$-parallel dividing curves 
which split off half-disks, where the half-disk along $\gamma\times \{1\}$ is 
positive (resp.\ negative) and there are $n$ parallel homotopically essential 
closed curves. See Figure~\ref{fig3} for the possibilities. 

\subsection{Type $\ii_n^\pm$, $n$ even.}  \label{section21}
The first case we treat is $\Gamma_A=\ii_n^\pm$, where $n$ is even.  

\begin{lemma}     \label{reduction1}
$\ii_{2m}^\pm$, $m\in \Z^+$, can be reduced to $\ii_0^\pm$.
\end{lemma}

\begin{proof}
The proof strategy is to start with the convex annulus $A$ with dividing set 
$\Gamma_A$ and then search for a bypass $B$ attached along $A$ such that 
isotoping $A$ across $B$ will produce a convex annulus $A'$ with $\Gamma_A'$ of 
decreased complexity.  It is important to keep in mind that Theorem 
\ref{basecase} is false for tori; thus in searching for $B$ we must exploit the 
assumption that the genus of $\Sigma$ is  $\geq 2$.

Assume we have $\ii_{2m}^+$. Figure~\ref{fig4} depicts the  
convex decomposition sequence for this case.  We will treat the case $m=1$, 
which is the hardest case.  The situation $m>1$ will be 
left to the reader.     

Figure~\ref{fig4}(A) depicts $M$ cut open along $A$.  Figure~\ref{fig4}(A) shows 
$(\Sigma\setminus \gamma)\times\{1\}$ together with $A^+$ and $A^-$, two copies 
of $A$ on $M\setminus A$.  (Warning:  $A^+$ and $A_+$ are distinct surfaces. $A_+\subset 
A$ is reserved for the positive region of a convex annulus $A$, whereas $A^+$ is 
the copy of $A$ in $M\setminus A$ where the induced orientation from $A$ agrees 
with the the orientation induced from $\bdry (M\setminus A)$.)  Figure 
\ref{fig4}(B) depicts $(\Sigma\setminus \gamma)\times\{0\}$. 
\begin{figure} [ht]	
	{\epsfysize=7in\centerline{\epsfbox{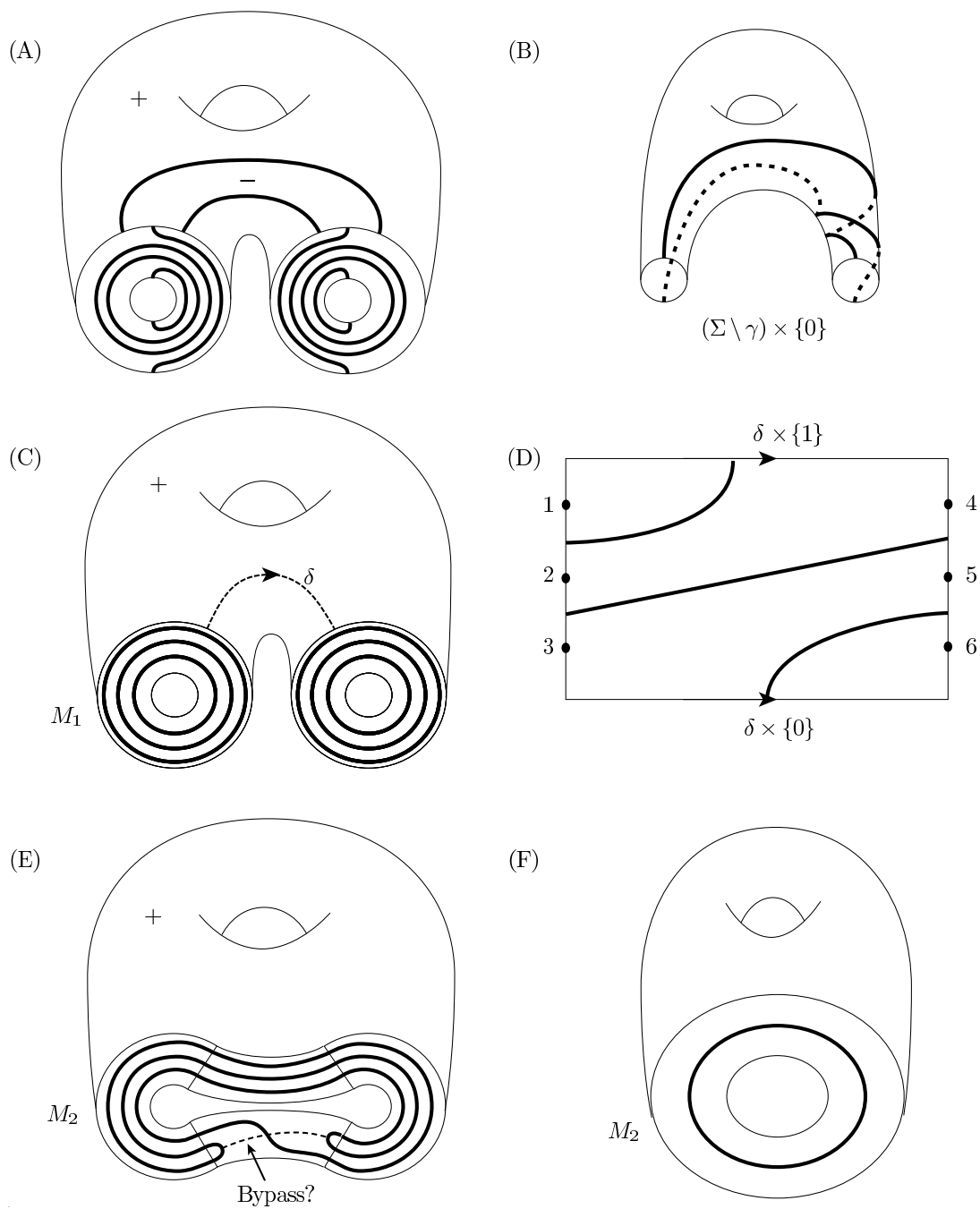}}}
	\caption{The case $\ii_{2m}^+$.} 	
	\label{fig4}
\end{figure}
After rounding the edges of $M\setminus A$, we obtain the convex handlebody 
$M_1$ in Figure~\ref{fig4}(C). The dividing set $\Gamma_{\bdry M_1}$ then 
consists of three parallel curves isotopic to the core curve of $A^+$ and three 
parallel curves isotopic to the core of $A^-$.

We now make the next cut in the convex decomposition along $\delta\times I$, 
where $\delta\subset \Sigma \setminus \gamma$ is a properly embedded oriented 
arc which connects the two boundary components of $\Sigma\setminus \gamma$ 
(from $A^+$ to $A^-$). (Some rounding will have taken place, but we assume that 
has already been taken care of.)   $\delta\times I$ will be given the 
orientation induced from $\delta$ and $I$.  We now consider the dividing curve 
configurations on $\delta\times I$.    $\bdry (\delta\times I)$ will intersect 
$\Gamma_{\bdry M_1}$ in three points along $A^+$ and in three points along 
$A^-$.  We label them 1, 2, 3 on $A^+$ in order from closest to $\delta\times 
\{1\}$ to farthest from $\delta\times \{1\}$.  Similarly label the three points 
of intersection on $A^-$ by 4, 5, 6 from closest to $\delta\times \{1\}$ to 
farthest  (see Figure~\ref{fig4}(D)).   We claim that if there exists a 
$\bdry$-parallel dividing curve straddling one of Positions 2, 3, 4, or 5, then 
the bypass corresponding to any of these positions, when considered back on $A$, 
would give a bypass along $A$ (from one of the sides) and a new convex annulus 
isotopic to $A$ with fewer dividing curves.  Positions 2 and 5 give rise to 
bypasses whose Legendrian arcs of attachment are contained in $A$ and which 
intersect three distinct curves of $\Gamma_A$.   A bypass at Positions 3 or 4, 
when traced back to Figure~\ref{fig4}(A), also yields a bypass along $A$ which 
reduces $\#\Gamma_A$.   To realize this, we apply Bypass Sliding 
(Lemma~\ref{bypass-sliding}).  Now, Figure~\ref{fig4}(D) is the only remaining 
dividing curve configuration for $\delta\times I$.

Therefore, we can either reduce from $\ii_{2}^+$ to $\ii_{0}^+$ or obtain the 
dividing set as in Figure~\ref{fig4}(D).  In the latter situation, we proceed by 
rounding the edges of $M_1\setminus (\delta\times I)$ to get $M_2$, 
depicted in Figure~\ref{fig4}(E). This, after the unique dividing curve is 
straightened, is equivalent to Figure~\ref{fig4}(F). In other words, 
$\Gamma_{\bdry M_2}$ consists of one curve, and it separates $\bdry M_2$.

We claim there exists a bypass from the interior of $\bdry M_2$ along  
$(\delta\times I)^-$ as depicted in Figure~\ref{fig4}(E).  This follows from 
using Lemma~\ref{lemma:super-lerp} with $F=\bdry M_2$. 
Once  we have the bypass, adding it to the exterior of $(\delta 
\times I)^+$ as shown in Figure~\ref{fig4}(E) forces the existence of 
bypasses in Positions~3 and 4.  Therefore, we can always reduce from 
$\ii^\pm_{2m}$ to $\ii^\pm_0$. 
\end{proof}

\begin{lemma}  \label{lemma:ii0}
$\Gamma_A=\ii_0^\pm$ extends uniquely to a tight contact structure on $M$.  It is 
universally tight.  
\end{lemma}

\begin{proof} 
After cutting $M$ along $A$, we obtain $M\setminus A=S\times I$, where $S$ is a 
surface of genus $g-1$ with two punctures.  Applying edge-rounding, we obtain 
that $\Gamma_{\bdry(S\times I)}$ is isotopic to $(\bdry S)\times \{{1\over 
2}\}$.  Lemma~\ref{uniqueness-bdry-case} (or the proof of 
Lemma~\ref{uniqueness-bdry-case}) implies that there is a unique tight contact 
structure which extends to the interior of $S\times I$.  The tight contact 
structure is universally tight by Lemma~\ref{uniqueness-bdry-case}, and glues to 
give a universally tight contact structure on $M$, since $\Gamma_A$ is 
$\bdry$-parallel and we can therefore apply Theorem~\ref{colin-gluing}. 
\end{proof}

\subsection{Type $\ii_n^\pm$, $n$ odd.}  \label{section22}

\begin{lemma}
$\ii_{2m+1}^\pm$, $m\in \Z^+$, can be reduced to $\ii^\pm_1$.
\end{lemma} 

\begin{figure} [ht]	
	{\epsfysize=4.5in\centerline{\epsfbox{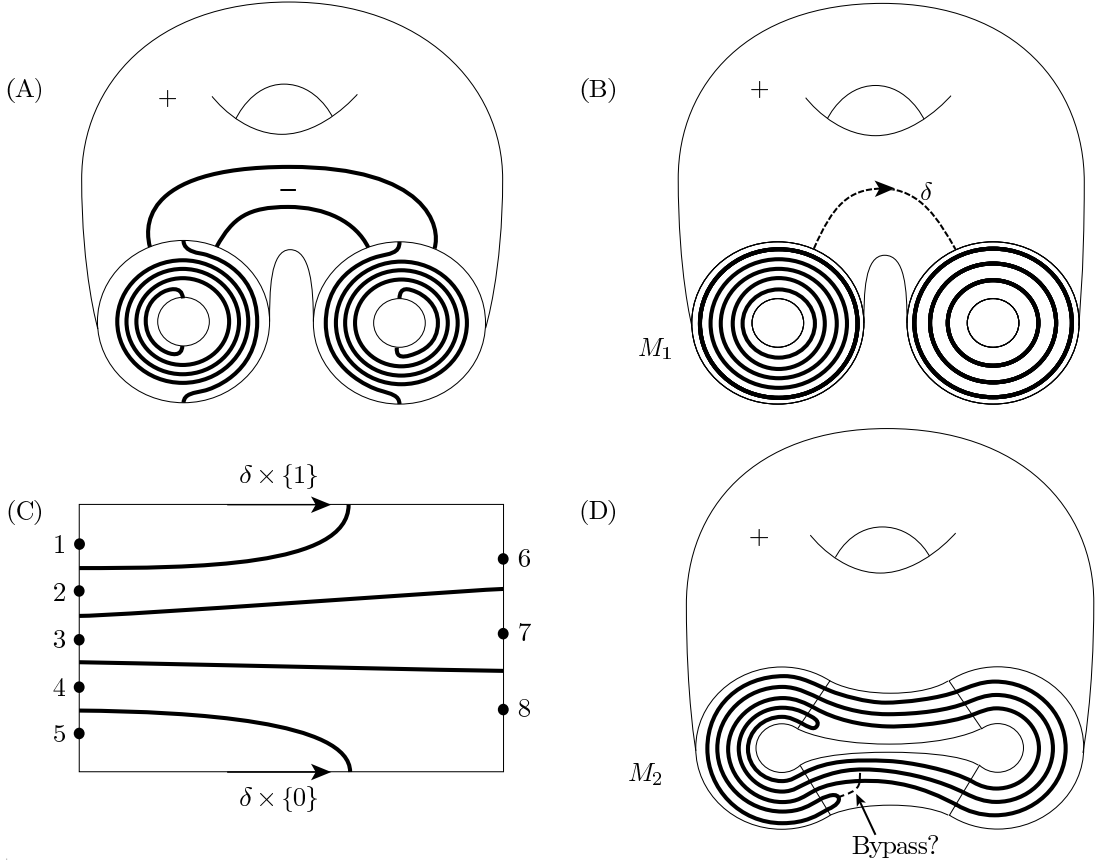}}}
	\caption{The case $\ii_{2m+1}^+$.} 	
	\label{fig5}
\end{figure}

\begin{proof}
The proof is similar to the proof of Lemma~\ref{reduction1}. It is enough to 
show that $\ii_{2m+1}^\pm$ can be reduced to $\ii_{2m-1}^\pm$.    Figure 
\ref{fig5}(A) depicts $M\setminus A$ where we have $\ii_{2m+1}^+$.  After 
rounding the edges, we obtain $M_1$ which has $2m+3$ closed curves parallel to 
the core curve of $A^+$ and $2m+1$ closed curves parallel to the core curve of 
$A^-$.  We take the next convex decomposing disk $\delta\times I$ with efficient
Legendrian boundary.   $\bdry(\delta\times I)$ intersects $\Gamma_{\bdry M_1}$ 
in $2m+3$ points along $A^+$, labeled $1$ through $2m+3$ from closest to 
$\delta\times\{1\}$ to farthest from $\delta\times\{1\}$, and $2m+1$ points 
along $A^-$, labeled $2m+4$ through $4m+4$.  The only $\bdry$-parallel dividing 
curves on $\delta\times I$ which do not immediately lead to a bypass on $A$ are 
those straddling Positions $1$ and $2m+3$.  Therefore, we are left to consider a 
unique choice for $\Gamma_{\delta\times I}$, given in Figure~\ref{fig5}(C), 
i.e., exactly two $\bdry$-parallel arcs (along $1$ and $2m+3$), and all other 
dividing curves consecutively nested around them.  

In order to prove the reduction from $\ii_{2m+1}^\pm$ to $\ii^\pm_{2m-1}$, we 
show the existence of a nontrivial bypass along $(\delta\times I)^-$ from the 
interior of $M_2$.  Indeed, Lemma~\ref{lemma:super-lerp} guarantees that there 
is a bypass along an arc of attachment which intersects the $\bdry$-parallel 
component of $(\delta\times I)^-$ straddling Position $1$, as well as two 
consecutively nested dividing arcs. See Figure~\ref{fig5}(D).  This bypass, if 
viewed on $\delta\times I$ (Figure~\ref{fig5}(C)), produces $\bdry$-parallel 
curves across Positions~$2m+2$ and $4m+3$, giving us a reduction in the number 
of parallel curves on $\Gamma_A$. \end{proof}

\begin{lemma}  \label{ii1}
$\ii^\pm_1$ can be reduced to ${\mbox{\it I}}_k$.
\end{lemma}

\begin{proof}      
We will treat $\Gamma_A=\ii^-_1$, and leave $\ii^+_1$ to the reader.  The 
computation is similar in spirit to the previous computations, except that it is 
a bit more involved.  The goal is to find a bypass along $A^+$ from the interior 
of $M_1$ which straddles the three components of $\Gamma_{A^+}$.  This time, the 
holonomy of the bypass (how many times the bypass wraps around the core curve) 
is important.  In order to determine the existence of a bypass, we will 
successively cut $M_1$, leaving $A^+$ untouched, until we arrive at a solid 
torus whose boundary contains $A^+$.  On the solid torus, we can determine 
whether the bypass exists, by appealing to the classification of tight contact 
structures on solid tori \cite{Gi00a, H1}.

\begin{figure}[ht]	
	{\epsfysize=4.5in\centerline{\epsfbox{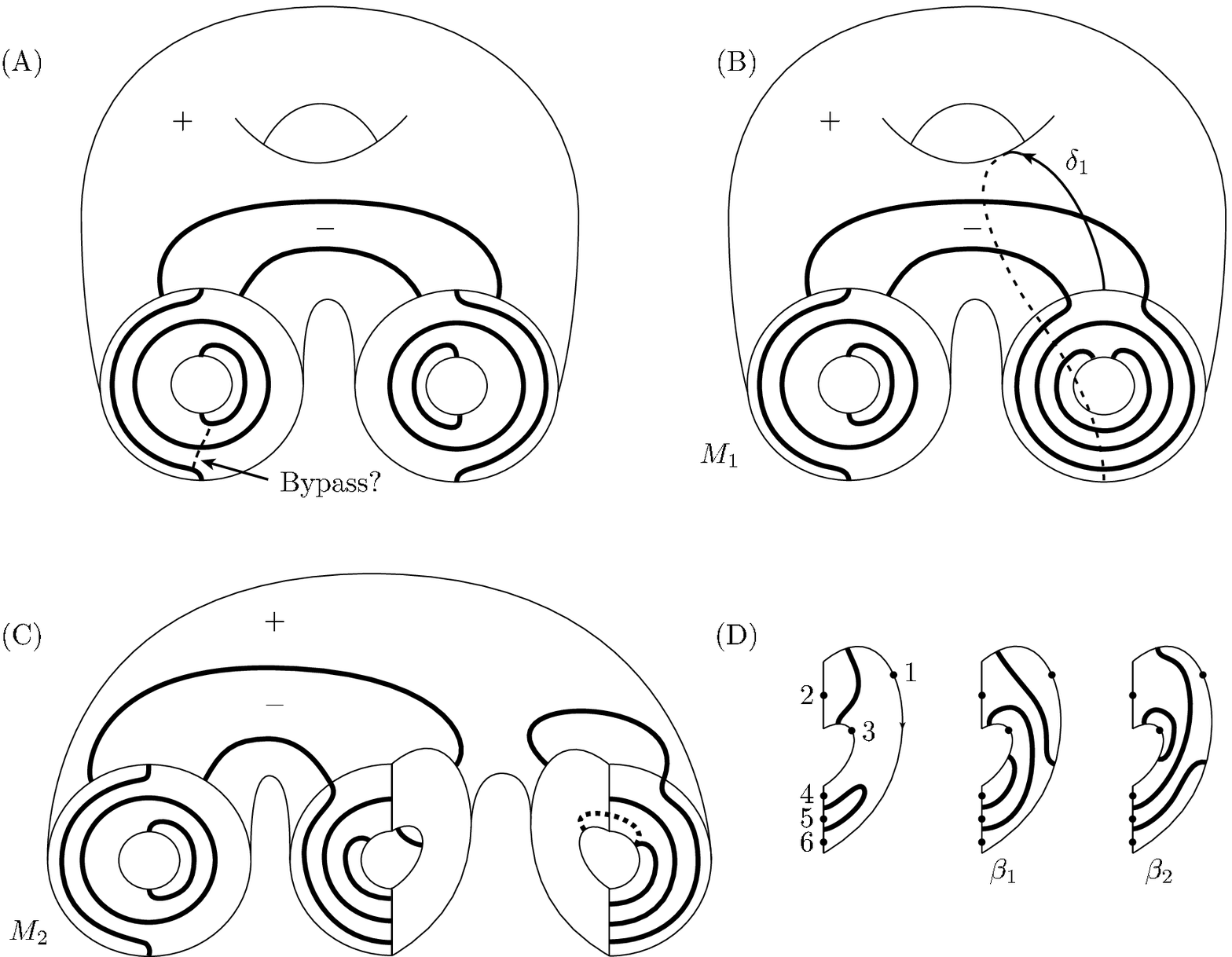}}}
	\caption{The case $\ii_1^-$.} 	
	\label{fig6}
\end{figure}

Figures~\ref{fig6} and \ref{fig7} depict this decomposition process.  As shown 
in Figure~\ref{fig6}(B), let $\delta_1$ be a non-boundary-parallel, properly 
embedded arc on $\Sigma\setminus \gamma$ which does not intersect $\delta$ and 
begins and ends on $\bdry A^-$.  We take $\bdry(\delta_1\times I)$ to be 
Legendrian and efficient with respect to $\Gamma_{\bdry M_1 \setminus A^+}$ on 
$\bdry M_1\setminus A^+$. Note this happens to be the same as being efficient 
with respect to $\Gamma_{\bdry M_1}$ on $\bdry M_1$.   The side of  Figure 
\ref{fig6}(A) which is hidden, namely $(\Sigma\setminus \gamma)\times \{0\}$, is 
as shown in Figure~\ref{fig4}(B), that is, the hidden side differs from the top 
by a single Dehn twist.  However, for subsequent figures, we suppose that the 
$\bdry$-parallel dividing curve on $A^-$ along $(\gamma\times \{0\})^-$ has 
already been twisted around the core curve of $A^-$, i.e., 
$\Gamma_{(\Sigma\setminus\gamma)\times\{0\}}$ is the same as 
$\Gamma_{(\Sigma\setminus\gamma)\times\{1\}}$. Now, take  $\delta_1\times I$ to 
be convex and define $M_2$ to be $M_1\setminus (\delta_1\times I)$, after 
rounding the edges.  (Warning: This $M_2$ is different from the $M_2$'s in the 
previous lemmas.)    We label  $\bdry (\delta_1\times I)\cap \Gamma_{\bdry M_1}$  
by numbers $1$ through $6$ as follows:  There is a unique closed curve of 
$\Gamma_{\bdry M_1\setminus A^+}$, and we let the two points of intersection of 
this curve with $\bdry (\delta_1\times I)$ be $2$ and $5$.  The rest are labeled 
in increasing order along $\bdry(\delta_1\times I)$ in the direction given by 
the induced orientation, i.e., ``counterclockwise''.  Now,  $\bdry$-parallel 
dividing curves in Positions $2$ and $5$ would immediately allow us to reduce to 
$I_k$, as can be seen on the $(\delta_1\times I)^-$-side.  Therefore, we are 
left with two possibilities for $\Gamma_{\delta_1\times I}$, which we call 
$\beta_1$ and $\beta_2$.  In Figure~\ref{fig6}(D), the left diagram represents 
the two $\bdry$-parallel positions which are ruled out, and the middle and right 
respectively are $\beta_1$ and $\beta_2$.

\begin{figure}[ht]	
	{\epsfysize=7in\centerline{\epsfbox{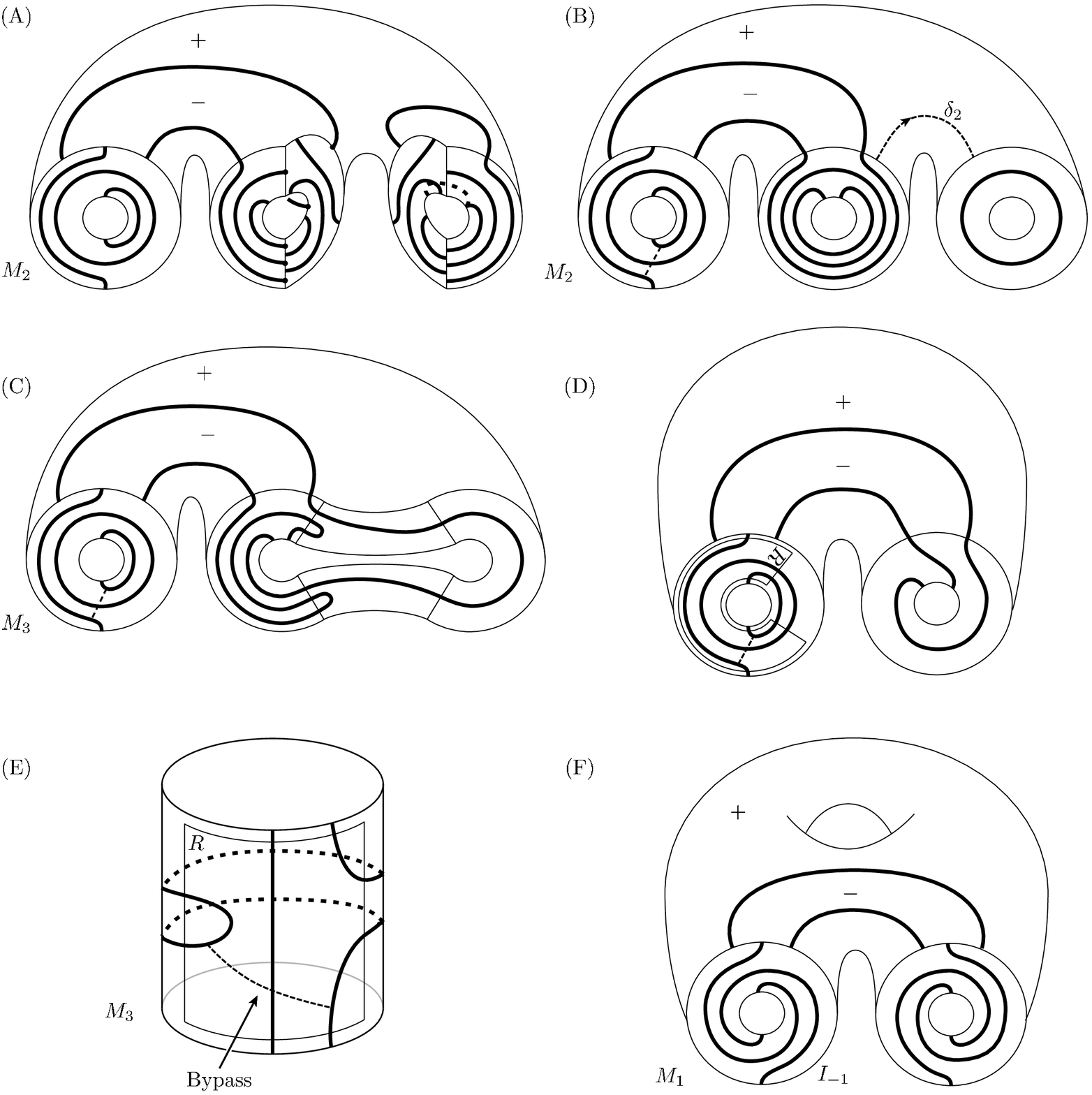}}}
	\caption{The case $\ii_1^-$, continued.  Configuration $\beta_1$.} 	
	\label{fig7}
\end{figure}

Consider $\beta_1$. (See Figure~\ref{fig7}.) Figure~\ref{fig7}(A) represents the 
dividing set $\Gamma_{\bdry M_2}$ before rounding, and Figure~\ref{fig7}(B) is 
$\Gamma_{\bdry M_2}$ after rounding.  We take $\delta_2$ (as in Figure~
\ref{fig7}(B)) to be a non-boundary-parallel, properly embedded arc on 
$\Sigma\setminus (\gamma\cup \delta_1)$ which begins on one copy of $\delta_1$ 
and ends on the other copy.  At this point, $\Sigma\setminus 
(\gamma\cup\delta_1)$ is a pair-of-pants, and cutting along $\delta_2$ yields an 
annulus.  Take $\bdry (\delta_2\times I)$ to be Legendrian and efficient with 
respect to $\Gamma_{\bdry M_2\setminus A^+}$  (which also happens to be the same 
as being efficient with respect to $\Gamma_{\bdry M_2}$).  Now, $tb(\bdry 
(\delta_2\times I))=-2$, and there are two possibilities for 
$\Gamma_{\delta_2\times I}$.  We may rule out one of the possibilities, since it 
yields a $\bdry$-parallel dividing curve which allows us to reduce to $I_k$.  

Finally, we have a solid torus $M_3=M_2\setminus (\delta_2\times I)$, and 
$\Gamma_{\bdry M_3}$ consists of 2 parallel longitudinal dividing curves. 
(Figures~\ref{fig7}(C,D,E) represent $\Gamma_{\bdry M_3}$ before and after 
successive edge-roundings.)  A compressing disk $D$ intersecting each dividing 
curve once may be chosen so that $M_3 \setminus D$ is as shown in 
Figure~\ref{fig7}(E), and in particular so that the rectangles labelled $R$ in 
Figure~\ref{fig7}(D) and Figure~\ref{fig7}(E) correspond.  This implies, by 
\cite{H1}, that $M_3$ is a standard neighborhood of a Legendrian curve.  Let us 
identify $\bdry M_3=\R^2/\Z^2$ by letting the meridian have slope $0$ and 
$\Gamma_{\bdry M_3}$ have slope $\infty$. 

We claim that there exists a bypass along $A^+$ from the interior of 
$M_1$, which changes $\ii^-_1$ to ${\mbox{\it I}}_{-1}$.  We look for the corresponding 
bypass along $A^+$ on $M_3$.  To see this exists, let $D$ be a convex meridional 
disk for $M_3$ with Legendrian boundary $\bdry D$ which is efficient with 
respect to $\Gamma_{\bdry M_3}$ and is disjoint from the bypass arc of 
attachment.  Then $tb(\bdry D)=-1$, and there is a unique way, up to isotopy, to 
cut this solid torus into a 3-ball $B^3$.  Finally, the bypass along $\bdry B^3$ 
is a trivial bypass, which must therefore exist by Right-to-Life.    
\end{proof}

\subsection{Type ${\mbox{\it I}}_k$.}    \label{section23}

\begin{lemma}     \label{lemma:ik1}
${\mbox{\it I}}_k$, $k> 0$, can be reduced to ${\mbox{\it I}}_0$.
\end{lemma}

\begin{proof}
The decomposition procedure is given in Figure~\ref{fig8}.  Figure~\ref{fig8}(A) 
gives $M\setminus A$;  note that in this figure $\Gamma_{(\Sigma\setminus 
\gamma)\times\{0\}}$ does not equal $\Gamma_{(\Sigma\setminus 
\gamma)\times\{1\}}$ and rather is as in Figure~\ref{fig4}(B).  
Figure~\ref{fig8}(B) is the same as Figure~\ref{fig8}(A), except that the extra 
Dehn twist is incorporated in $A^-$ so that $\Gamma_{(\Sigma\setminus 
\gamma)\times\{0\}}=\Gamma_{(\Sigma\setminus \gamma)\times\{1\}}$.   
Figure~\ref{fig8}(C) depicts $\bdry M_1$, where $M_1$ is $M\setminus A$ with the 
edges rounded.  

\begin{figure}[ht]	
	{\epsfysize=4in\centerline{\epsfbox{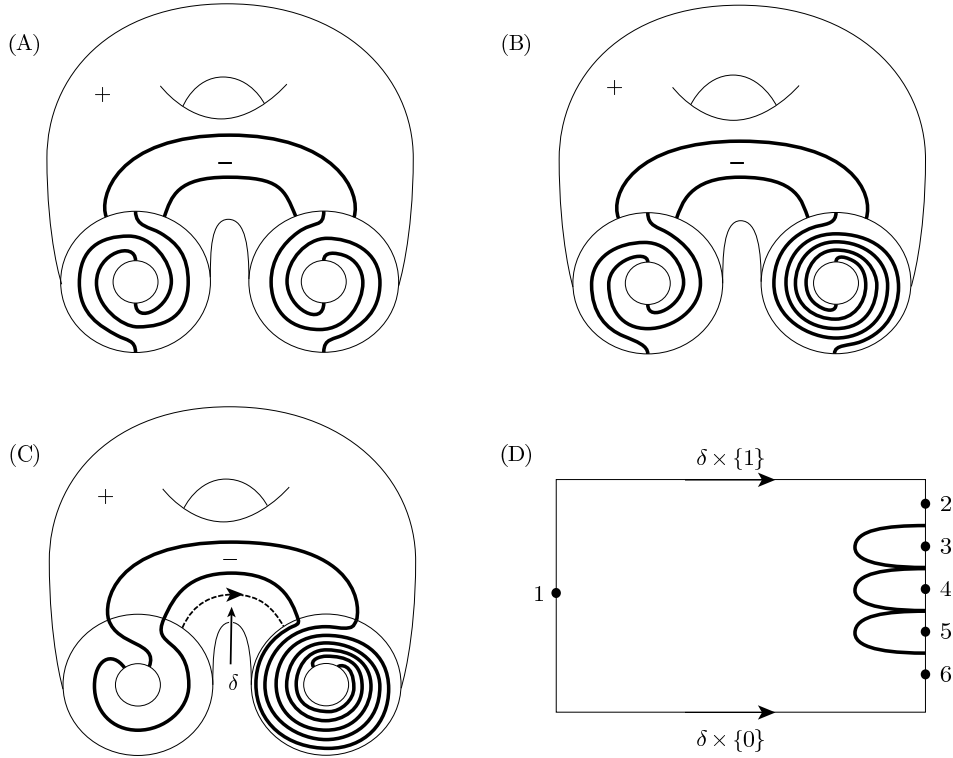}}}
	\caption{The case ${\mbox{\it I}}_k$, $k> 0$.} 	
	\label{fig8}
\end{figure}

Let $\delta\subset \Sigma\setminus\gamma$ be the same as in Lemma~
\ref{reduction1}.  As before, consider the compressing disk $\delta\times I$, 
which we take to be convex with efficient Legendrian boundary.  $\delta\times 
I$ is chosen so that $\bdry (\delta\times I)$ intersects $\Gamma_{\bdry M_1}$ 
in $2k-1$ points along $A^+$ (labeled $1$ through $2k-1$ in order from closest 
to $\delta\times \{1\}$ to farthest) and $2k+3$ points along $A^-$ (labeled  
$2k$ through $4k+2$ in order from closest to $\delta\times\{1\}$ to farthest).  
If there are $\bdry$-parallel components of $\Gamma_{\delta\times I}$ along any 
of $2k+1,\dots,4k+1$, then the corresponding bypasses would give rise to the 
state transition from ${\mbox{\it I}}_k$ to ${\mbox{\it I}}_{k-1}$.  Now, there are $2k+2$ endpoints of 
$\Gamma_{\delta\times I}\cap \bdry (\delta\times I)$ between Positions $2k$ and 
$4k+2$.   If there are connections (dividing arcs) amongst the  $2k+2$ 
endpoints, then clearly, this would give rise to a $\bdry$-parallel arc 
straddling one of the ``reducing'' positions.  However, this must happen since 
the total number of endpoints is $4k+2$, i.e., $\#(\Gamma_{\delta\times I}\cap 
\bdry (\delta\times I))=4k+2$, and $4k+2<2(2k+2)$. \end{proof}

\begin{lemma}     \label{lemma:ik2}
${\mbox{\it I}}_k$, $k<-1$, can be reduced to ${\mbox{\it I}}_{-1}$.
\end{lemma}

The apparent lack of symmetry between Lemmas~\ref{lemma:ik1} and \ref{lemma:ik2}
is due to the fact that the first cut along $A$ is not symmetric with respect to 
$\Gamma_0$.

\begin{figure}[ht]	
	{\epsfysize=4in\centerline{\epsfbox{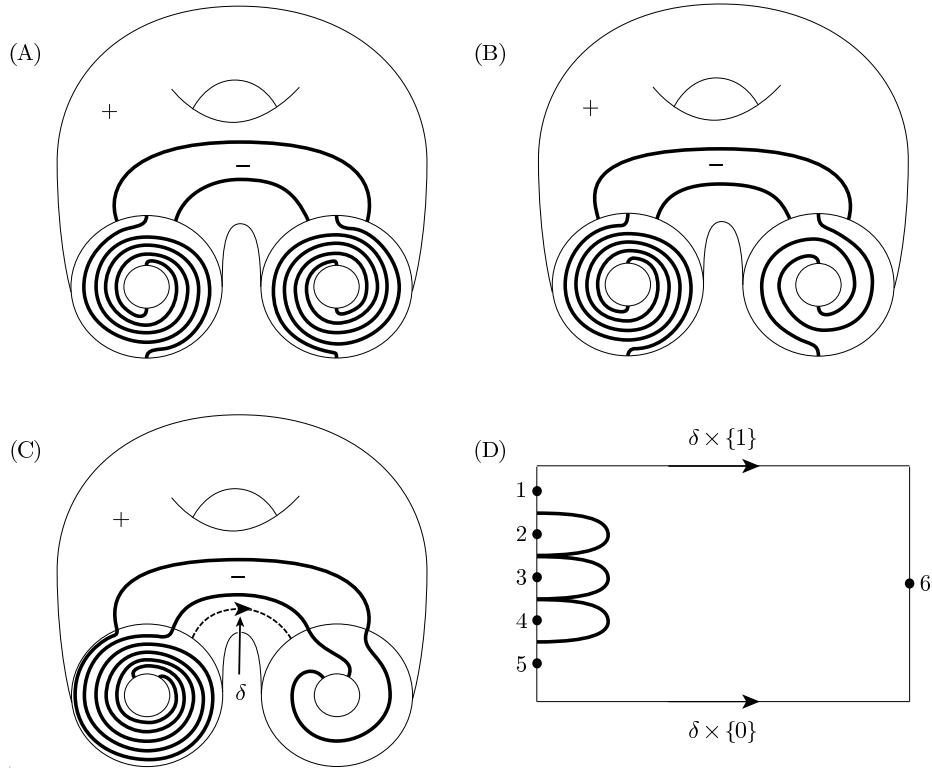}}}
	\caption{The case ${\mbox{\it I}}_k$, $k<-1$.} 	
	\label{fig9}
\end{figure}

\begin{proof} 
The argument is almost identical to that of Lemma~\ref{lemma:ik1}.  The only 
difference is that the computations are mirror images of those of Lemma~
\ref{lemma:ik1}.  Refer to Figure~\ref{fig9} for the steps of the computation.
\end{proof}

Lemmas~\ref{lemma:ik1} and \ref{lemma:ik2} together indicate that all the 
${\mbox{\it I}}_k$'s can be reduced to either ${\mbox{\it I}}_0$ or ${\mbox{\it I}}_{-1}$.  The following Lemma~
computes (an upper bound for) the tight contact structures where 
$\Gamma_A={\mbox{\it I}}_0$, and relates them to ${\mbox{\it I}}_{-1}$.

\begin{lemma}
There are at most two tight contact structures on $\Sigma\times I$ in the Base 
Case for which $\Gamma_A={\mbox{\it I}}_0$.  The same also holds for $\Gamma_A={\mbox{\it I}}_{-1}$.  
There exist state transitions along $A$ which allow us to switch between ${\mbox{\it I}}_0$ 
and ${\mbox{\it I}}_{-1}$.
\end{lemma}

\begin{proof}
Let us take the case $\Gamma_A={\mbox{\it I}}_0$. Figure~\ref{fig10}(A) gives the dividing 
set of $M\setminus A$, before edge-rounding but after the extra Dehn twist from 
the bottom face $(\Sigma\setminus \gamma)\times\{0\}$ is included. (Therefore,
$\Gamma_{(\Sigma\setminus\gamma)\times\{0\}} 
=\Gamma_{(\Sigma\setminus\gamma)\times\{1\}}$.)         Figure~\ref{fig10}(B) is 
the same after edge-rounding.  Now, if we cut along $\delta\times I$, defined as 
in Lemma~\ref{lemma:ik1}, there are two possibilities for $\Gamma_{\delta\times 
I}$, since $tb(\bdry(\delta\times I))=-2$. (These are shown in Figures 
\ref{fig10}(C,D).)  Figures~\ref{fig10}(E,F) depict the dividing set of 
$M_2=M_1\setminus (\delta\times I)$, after edge-rounding.  In both cases, 
$\Gamma_{\bdry M_2}$ consists of exactly one dividing curve parallel to $\bdry 
(\Sigma\setminus (\gamma\cup \delta))$.   Finally, using Lemma~
\ref{uniqueness-bdry-case}, we find that for each of the two possibilities of 
$\Gamma_{\delta\times I}$ there is a unique universally tight contact structure 
on $M_2$.  Theorem~\ref{colin-gluing} is also sufficient to glue back along 
$\delta\times I$ to give two universally tight contact structures on $M_1$ with 
the boundary condition given by Figure~\ref{fig10}(A). Making the final gluing 
and proving the resulting contact structure is universally tight is a more 
complicated state-transition operation, so we will content ourselves for the 
time being with the knowledge that there are at most two tight contact 
structures on $\Sigma\times I$ in the Base Case with $\Gamma_A={\mbox{\it I}}_0$.      A 
similar computation also holds for $\Gamma_A={\mbox{\it I}}_{-1}$.  

\begin{figure}[ht]	
	{\epsfysize=7in\centerline{\epsfbox{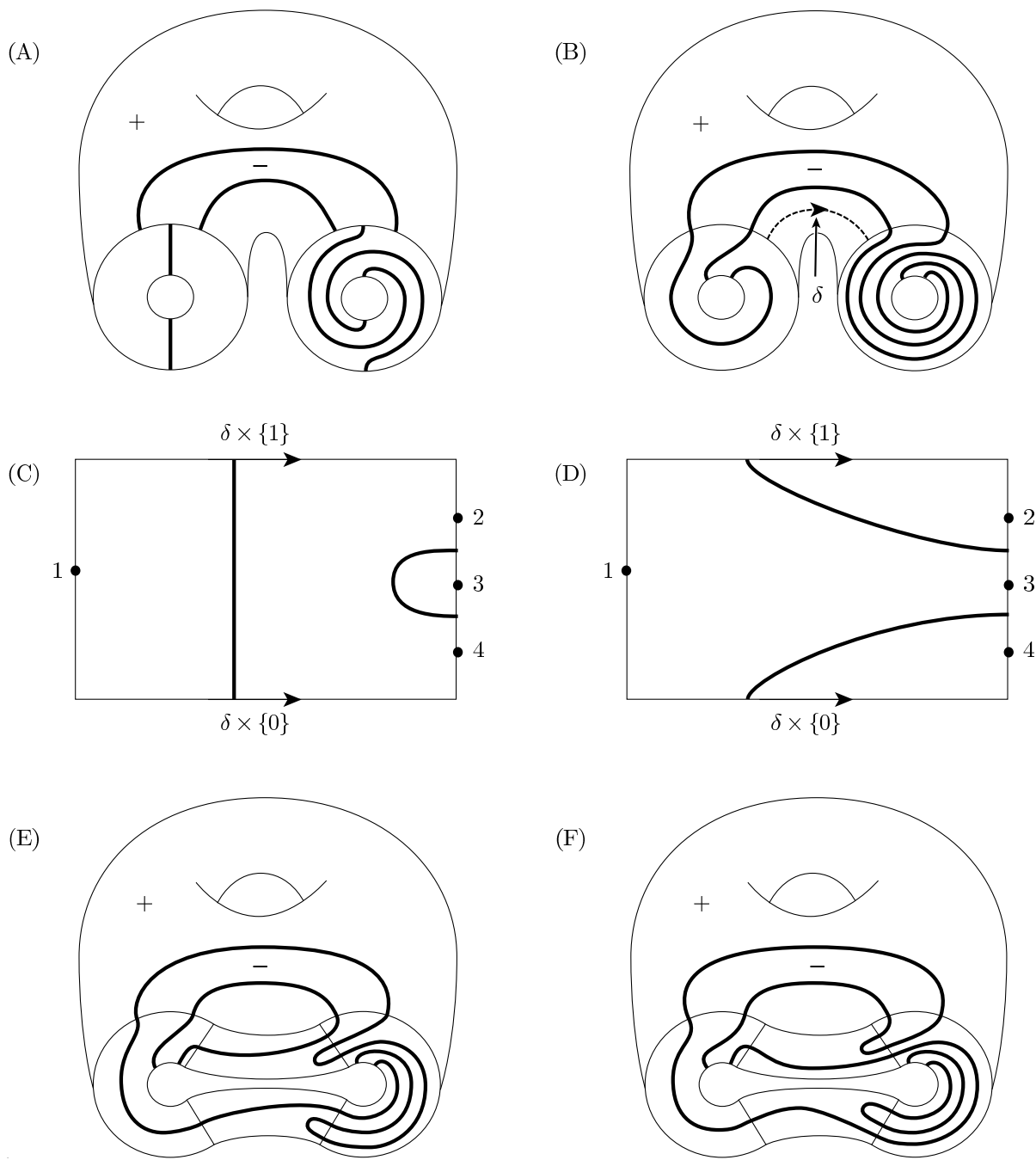}}}
	\caption{Classification for ${\mbox{\it I}}_0$.} 	
	\label{fig10}
\end{figure}

We now prove that we may switch from ${\mbox{\it I}}_0$ to ${\mbox{\it I}}_{-1}$.  The reverse 
procedure is identical.  Above, we found that $\#(\Gamma_{\bdry M_1}\cap \bdry 
(\delta \times I))=4$, and one intersection was on $A^+$ (labeled $1$) and $3$ 
on $A^-$ (labeled $2,3,4$ in succession from closest to $\delta\times\{1\}$ to 
farthest). A $\bdry$-parallel dividing curve of $\delta\times I$ straddling 
Position $3$ clearly allows us to transition from ${\mbox{\it I}}_0$ to ${\mbox{\it I}}_{-1}$.  On the 
other hand, a $\bdry$-parallel dividing curve straddling Position $4$ gives rise 
to a corresponding bypass which can be slid so that all three of the 
intersections with $\Gamma_{\bdry M_1}$ lie on $A^-$.   Therefore, for both 
choices of $\Gamma_{\delta\times I}$, we may transition from ${\mbox{\it I}}_0$ to 
${\mbox{\it I}}_{-1}$.
\end{proof}

\subsection{Completion of Theorem~\ref{basecase}.}   \label{section24}
Summarizing what we have proved so far:

\s\n
\begin{itemize}
\item $\#\pi_0(Tight(M,\mathcal{F}))\leq 4.$
\item Type $\ii^+_0$: there exists one.
\item Type $\ii^-_0$: there exists one.  
\item Type ${\mbox{\it I}}_0$: there are at most $2$.
\item The other possibilities for $\Gamma_A$ reduce to one of $\ii^\pm_0$ or 
${\mbox{\it I}}_0$. 
\item The tight contact structures of type $\ii^\pm_0$ are universally 
tight by Lemma~\ref{lemma:ii0}.  
\item The tight contact structures of type $\ii^+_0$ and $\ii^-_0$ are distinct 
and are also distinct from those of type ${\mbox{\it I}}_0$ by the relative Euler class 
$\tilde e$ evaluated on $A$. 
\end{itemize}

\s\n
The proof of Theorem~\ref{basecase} is complete, once we show Lemmas~\ref{first} 
and \ref{second} below. 

\begin{lemma}\label{first}
There are exactly two tight contact structures of type ${\mbox{\it I}}_0$.  They are 
universally tight, are obtained by adding a single bypass onto $\Sigma_0$, 
and have relative Euler class $PD(\tilde e(\xi))=\pm(\gamma_1-\gamma_0)$.
\end{lemma}

\begin{figure}[ht]	
	{\epsfysize=5in\centerline{\epsfbox{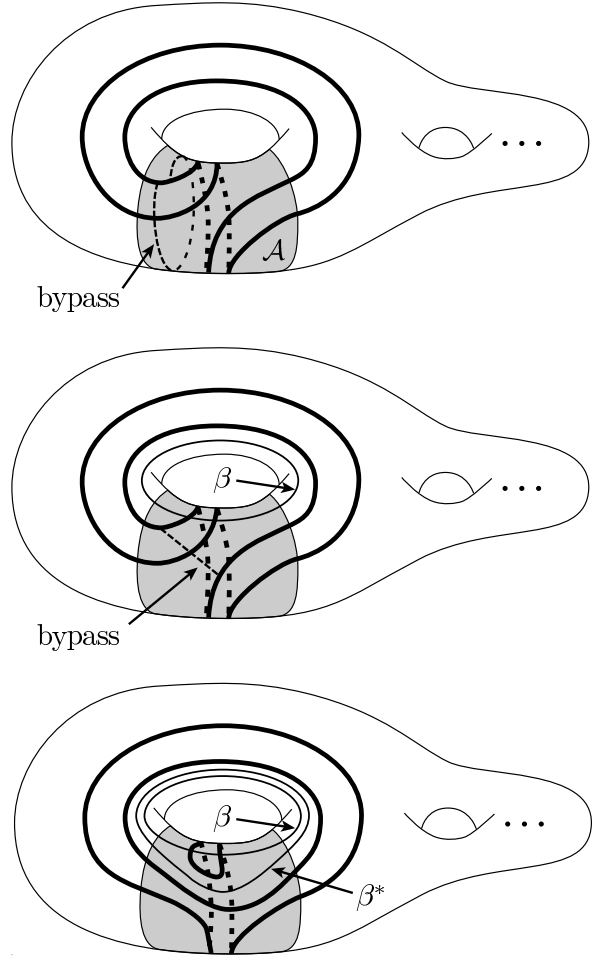}}}
	\caption{Ocarinian Riemann surfaces} 	
	\label{fig13}
\end{figure}

\begin{proof}  We show the existence of the (universally) tight contact 
structures by embedding them inside a suitable (universally) tight contact 
structure of type $\ii^\pm_0$.  Let $A=\gamma\times I$ be the first cut that was 
used to decompose $\Sigma\times I$. If $\Gamma_A=\ii^+_0$, then there is a 
$\bdry$-parallel dividing curve along $\gamma\times \{0\}$ which cuts off a 
positive region of $A$.  There is a corresponding degenerate bypass whose 
Legendrian arc of attachment is all of $\gamma\times\{0\}$.  (Degenerate means 
that the two ends of the Legendrian arc of attachment are identical.)  If we 
immediately attach this degenerate bypass onto $\Sigma_0$, we obtain an isotopic 
convex surface which we call $\Sigma_{1/2}'$ and which satisfies 
$\Gamma_{\Sigma_{1/2}'}=2\gamma$.   However, if we separate the endpoints by 
bypass sliding in one particular direction along $\Gamma_{\Sigma_0}$, we 
obtain a more convenient convex surface $\Sigma_{1/2}$ with 
$\Gamma_{\Sigma_{1/2}}=2\gamma_1$, after the bypass attachment.  There are two 
possible bypasses (of opposite sign) we can attach to $\Sigma_0$, arising from 
$\ii^+_0$ and $\ii^-_0$.  They are clearly universally tight by construction.

We now compute their relative Euler class.  It will then be clear that the two 
tight contact structures are distinct and of type ${\mbox{\it I}}_0$.   We fix some 
notation.  The tight contact structure $\xi$ obtained in the previous paragraph 
by attaching a bypass has ambient manifold $\Sigma\times [0,1]$, 
$\Gamma_{\Sigma_i}=2\gamma_i$, $i=0,1$, and $\xi$ is $[0,1]$-invariant except 
for $\mathcal{A}\times [0,1]$, where $\mathcal{A}\subset \Sigma_0$ is a convex 
annulus with Legendrian boundary whose core curve is isotopic to $\gamma$, and 
the bypass was attached to $\Sigma_0$ along $\mathcal{A}$.  

First suppose $\beta\subset \Sigma$ is a closed nonseparating curve which 
intersects neither $\gamma_0$ nor $\gamma_1$.  Then the corresponding convex 
annulus $\beta\times I$ will only consist of closed curves parallel to the core 
curve.  Hence $\langle \tilde{e}(\xi),\beta\times I\rangle =0$. 
Thus we may choose a basis of $H_1(\Sigma;\Z)$ such that, of the $2g$ 
generators,  $2g-2$ of them evaluate to zero in this manner.  Next let 
$\beta\times \{0\}\subset \Sigma_0$ be a closed Legendrian curve parallel to 
$\gamma$ but disjoint from $\mathcal{A}$.  Then, since $\xi$ is $I$-invariant 
away from $\mathcal{A}\times I$, $\beta\times I$ must consist of $2$ parallel 
vertical nonseparating arcs, that is, $\langle \tilde{e}(\xi),\beta\times 
I\rangle =0$.  Finally, let $\beta\times \{0\}\subset \Sigma_0$ be a closed 
efficient Legendrian curve which is parallel to $\gamma_1$ but does not 
intersect the arc of attachment of the bypass, which we take to be 
nondegenerate. Then $\beta\times\{i\}$, $i=0,1$, intersects $\Gamma_{\Sigma_i}$ 
twice,  and $\Gamma_{\beta\times I}$ consists of 2 parallel vertical 
nonseparating arcs.  However, now $\beta\times \{1\}$ is not efficient with 
respect to $\Gamma_{\Sigma_i}$, and resolving the extra intersection to produce 
an efficient intersection gives $\langle \tilde{e}(\xi),(\beta\times 
I)^*\rangle=\pm 1$.  (Here, $(\cdot)^*$ refers to the annulus with efficient 
Legendrian boundary.)  See Figure~\ref{fig13}.  Having evaluated 
$\tilde{e}(\xi)$ on all the basis elements, we find that $PD(\tilde{e}(\xi))= 
\pm \gamma=\pm (\gamma_1-\gamma_0)$.  \end{proof}

\begin{defn}
A tight contact structure on $\Sigma\times I$ which is contact diffeomorphic to 
one of the tight contact structures of type ${\mbox{\it I}}_0$ is said to be a {\em basic 
slice}.  Thus, in a basic slice, $\Gamma_{\Sigma\times\{0\}}$ and 
$\Gamma_{\Sigma\times\{1\}}$  each consist of two parallel curves 
$2\gamma_0$ and $2\gamma_1$ with  $|\gamma_0\cap \gamma_1|=1$.  The 
basic slice is obtained by attaching a single bypass  $B$ onto 
$\Sigma\times \{0\}$ and thickening $(\Sigma\times \{0\})\cup B$.
\end{defn}

\begin{lemma}\label{second}
The tight contact structure of type $\ii^\pm_0$ has relative Euler class
$PD(\tilde e(\xi))=\mp (\gamma_1+\gamma_0)$.
\end{lemma}

\begin{proof}
We will restrict our attention to $\ii^+_0$.  
As in the proof of Lemma~\ref{first}, if $\beta\subset \Sigma$ is a closed 
nonseparating curve which does not intersect $\gamma_0$ or $\gamma_1$, then 
$\langle \tilde{e}(\xi),\beta\times I\rangle =0$.  From the definition of  
$\ii^+_0$ we have: 
$$\langle \tilde{e}(\xi),(\gamma_1-\gamma_0)\times 
I\rangle=2.$$ 
Next, since $\gamma_0\times I$ with efficient Legendrian boundary 
intersects $\Gamma_{\Sigma_0}$ $0$ times and $\Gamma_{\Sigma_1}$ twice, 
$$\langle \tilde{e}(\xi),\gamma_0\times I\rangle=\pm 1,$$ 
where the exact sign will be determined in a moment. 
Similarly,
$$\langle \tilde{e}(\xi),\gamma_1\times I\rangle=\pm 1.$$  
For the three equations to agree, we must have $\langle 
\tilde{e}(\xi),\gamma_0\times I\rangle=-1$ and $\langle 
\tilde{e}(\xi),\gamma_1\times I\rangle=1$.  This implies that  $PD(\tilde 
e(\xi))=- (\gamma_1+\gamma_0)$. 
\end{proof}

\section{Classification of tight contact structures on $\Sigma\times I$}  
\label{section:classification}

In this section we prove Theorem~\ref{classification} as well as the following 
theorem:

\begin{thm}[Gluing Theorem]\label{gluing}
Let $\Sigma$ be an oriented closed surface of genus $g \geq 2$, $M=\Sigma\times 
[0,2]$, and $\xi$ a contact structure which is tight on $M\setminus \Sigma_1$.  
Suppose $\Sigma_i$, $i=0,1,2$, are convex and  $\Gamma_{\Sigma_i}=2\gamma_i$, 
where $\gamma_i$ are nonseparating oriented curves.   Also assume that 
the $\gamma_i$ are not mutually homologous.  If $PD(\tilde 
e(\xi|_{\Sigma\times[0,1]}))=\gamma_1-\gamma_0$ (here $\gamma_0$ and $\gamma_1$ 
have been oriented so the relative Euler class has this form), then $PD(\tilde 
e(\xi|_{\Sigma\times[1,2]}))=\gamma_2-\gamma_1$ (for some orientation of 
$\gamma_2$) if and only if $\xi$ is tight on $M$. \end{thm}

The condition of the $\gamma_i$ being mutually nonhomologous is a technical 
condition, which can be removed if we reformlate the Gluing Theorem without 
reference to the relative Euler class.  The reader is encouraged to do so, after 
examining the proof of Theorem~\ref{classification} and the Gluing Theorem.  As 
we will see, the only contact topology calculations needed to prove 
Theorem~\ref{classification} are the one done in Section~\ref{section:basecase} 
and a similar calculation in Proposition~\ref{exists-bypass}.  The rest 
is largely a ``proof by pure thought'', relying on the relative Euler class 
consistency check, Proposition~\ref{ubiquity}, and curve complex facts.

\subsection{Freedom of choice}
The proof of Theorem~\ref{classification} is founded on the following rather 
remarkable proposition.

\begin{prop} [Freedom of choice] \label{ubiquity}
Let $(M=\Sigma\times I,\xi)$ be a basic slice with $\Gamma_{\Sigma_i} = 
2\gamma_i$, $i=0,1$ and $\langle \gamma_1,\gamma_0\rangle=+1$, where $\langle 
\cdot \rangle$ represents the intersection form on $\Sigma$, a surface of genus 
at least $2$. Let $\gamma$ be any nonseparating curve on $\Sigma$.  Then there 
exists a convex surface isotopic to $\Sigma_0$ (or, equivalently, to 
$\Sigma_1$), which we call $\Sigma_{1/2}$ and which has 
$\Gamma_{\Sigma_{1/2}}=2\gamma$. 
\end{prop}

The strategy of the proof is to start with $\Gamma_{\Sigma_0}=2\gamma_0$ and
$\Gamma_{\Sigma_1}=2\gamma_1$ and successively find subslices $\Sigma\times 
[a,b]\subset \Sigma\times [0,1]$ with convex boundary and dividing sets 
$\Gamma_{\Sigma_a}$, $\Gamma_{\Sigma_b}$ which consist of two parallel 
curves each and represent curves which are ``closer'' to $\gamma$ inside the 
curve complex.     To prevent our notation from becoming to cumbersome, we will 
rename the old $\Sigma\times[a,b]$ to be the new $\Sigma\times[0,1]$ after each 
step of the induction.  We also write $[\alpha_0,\alpha_1; PD(\tilde{e}(\xi))]$ 
to mean some tight contact structure on $\Sigma\times [0,1]$ with  
$\Gamma_{\Sigma_i}=2\alpha_i$, $i=0,1$ and relative Euler class 
$PD(\tilde{e}(\xi))$.  If we do not specify the relative Euler class (or it is 
understood), we simply write $[\alpha_0,\alpha_1]$.  Moreover, if we want to 
indicate a basic slice, we write $\ll\alpha_0,\alpha_1\rr$.

We first describe the operation which will be used repeatedly in the proof.  

\begin{op}  \label{op}
Consider $[\alpha_0,\alpha_1]$, where $|\alpha_0\cap \alpha_1|=1$. Let $\alpha$ 
be a closed (necessarily nonseparating) curve which satisfies $|\alpha_0\cap 
\alpha|=0$ and $|\alpha_1\cap \alpha|=1$.  Then there exists a convex surface 
$\Sigma_{1/2}$ with $\Gamma_{\Sigma_{1/2}}=2\alpha$. 
\end{op}

\begin{proof}[Proof of Operation.]
On $\Sigma_0$, use LeRP to realize $\alpha\times\{0\}$ as a Legendrian curve 
with $\#(\Gamma_{\Sigma_0}\cap \alpha)=0$.  Similarly, Legendrian realize 
$\alpha\times\{1\}$ with $\#(\Gamma_{\Sigma_1}\cap \alpha)=2$.  Take the convex 
annulus $\alpha\times I$.  By the Imbalance Principle of \cite{H1}, there must 
be a $\bdry$-parallel dividing curve along $\alpha\times\{1\}$ and hence a 
degenerate bypass.  Attaching the degenerate bypass gives an isotopic convex 
surface with $\Gamma=2\alpha$.   
\end{proof}

\begin{proof}[Proof of Proposition~\ref{ubiquity}.]
We start with $\ll\gamma_0,\gamma_1\rr$.  Using Proposition~\ref{fact0}, we 
obtain a sequence:
$$\gamma_0=\alpha_0, \gamma_1=\alpha_1, \alpha_2,\dots, \gamma=\alpha_k,$$
where $\alpha_i$, $i=0,\dots,k$, are nonseparating, 
$|\alpha_{i-1}\cap\alpha_{i}|=1$, $i=1,\dots,k$, and 
$|\alpha_{i-1}\cap\alpha_{i+1}|=0$, $i=1,\dots,k-1$.  (Note that the statement 
of the proposition does not quite give what we want, but the proof clearly 
does.) Now, using the Operation, we successively find: 
$$\ll\alpha_0,\alpha_1\rr\supset \ll\alpha_2,\alpha_1\rr\supset 
\ll\alpha_2,\alpha_3\rr\supset \dots \supset \ll\alpha_{k-1},\alpha_k\rr.$$  
This completes the proof of the proposition.
\end{proof}

\subsection{Proof of Part 3 of Theorem~\ref{classification}.}  We will now give 
the classification result for $[\gamma_0,\gamma_1=\gamma_0]$, where $\gamma_0$ 
is an oriented  nonseparating curve.  If the tight contact structure $\xi$ is 
$I$-invariant, then $PD(\tilde{e}(\xi))=0$.

\begin{prop}  \label{exists-bypass}
Let $[\gamma_0,\gamma_0]$ be a tight contact structure on $M=\Sigma\times [0,1]$ 
which is not $I$-invariant.  Then $\Sigma\times [0,1]$ contains some basic slice 
$\Sigma\times[a,b]$.  
\end{prop}

\begin{proof}
The proof is a calculation along the lines of Section~\ref{section:basecase}.   
Let $\gamma\subset \Sigma$ be an oriented curve so that $\langle 
\gamma,\gamma_0\rangle=+1$.  Let $A=\gamma\times I$ be a convex annulus with 
efficient Legendrian boundary.  $\Gamma_A$ has the possibilities as in 
Figure~\ref{fig3}, denoted types ${\mbox{\it I}}_k$ and $\ii^\pm_n$.  Types $\ii^\pm_n$ all 
have $\bdry$-parallel dividing curves, which give rise to bypasses which, in 
turn, yield basic slices.  Therefore, it suffices to consider ${\mbox{\it I}}_k$.  

If $\Gamma_A={\mbox{\it I}}_0$, then there exists a sequence of convex meridional disks 
$D_i$ which decompose $M_1=M\setminus A$ into the 3-ball, and for which 
$tb(\bdry D_i)=-1$.  Hence, there must be a unique tight contact structure with 
$\Gamma_A={\mbox{\it I}}_0$.  This implies that $\Gamma_A={\mbox{\it I}}_0$ represents the $I$-invariant 
case.  

\begin{figure}[ht]	
	{\epsfysize=4in\centerline{\epsfbox{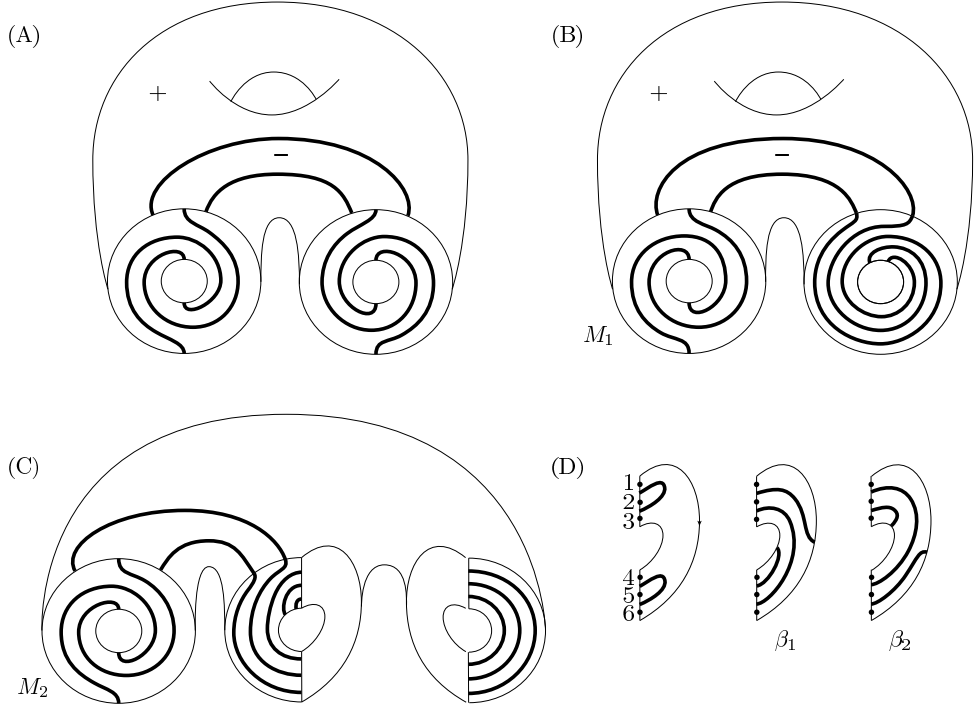}}}
	\caption{Transition from ${\mbox{\it I}}_k$ to $\ii^\pm_n$.} 	
	\label{fig11}
\end{figure}

Suppose $\Gamma_A={\mbox{\it I}}_k$ with $k>0$. Let $M_1=M\setminus A$.   
Figure~\ref{fig11}(A) gives $M_1$ and Figure~\ref{fig11}(B) the same with the 
edges of $A^-$ rounded.  (Note that Figure~\ref{fig11} depicts the case 
$\Gamma_A=I_1$.) They are presented in almost identical fashion as in 
Section~\ref{section:basecase} with one exception: now $\Gamma_{(\Sigma\setminus 
\gamma)\times \{1\}}=\Gamma_{(\Sigma\setminus \gamma)\times \{0\}}.$  Our 
notation will be identical to that of Lemma~ \ref{ii1}.  Let $\delta_1$ be an 
arc in $\Sigma \setminus \gamma$ with  $\bdry \delta_1= p_1-p_0$, where $p_0 \in 
A^+$, $p_1 \in A^-$, and $\delta_1$ does not intersect $\Gamma_{\Sigma}$. Take 
$\delta_1\times I$ and perturb it to be  convex with  Legendrian
boundary so that  $|\bdry (\delta_1\times I)\cap \Gamma_{\bdry M_1}|=4k+2$, and
$2k+1$ of the intersections we may assume are on $p_0\times I$ (labeled 
$1,\dots,2k+1$ from closest to $\delta_1\times\{1\}$ to farthest) and the other 
$2k+1$ on $p_1\times I$ (labeled $2k+2,\dots,4k+2$ from farthest from 
$\delta_1\times\{1\}$ to closest).  If there is a $\bdry$-parallel dividing arc 
on $\delta_1\times I$ straddling Positions $2,\dots,2k$ or $2k+3,\dots,4k+1$, 
then the corresponding bypass state transitions us into ${\mbox{\it I}}_{k-1}$.  If we can 
continue this, we eventually get to ${\mbox{\it I}}_0$, which is already taken care of. 
(Actually, it is unlikely such a state transition exists; we probably have an 
overtwisted contact structure here.)  Otherwise, we have two possibilities: 
$\beta_1$ and $\beta_2$.

\begin{figure}[ht]	
	{\epsfysize=4in\centerline{\epsfbox{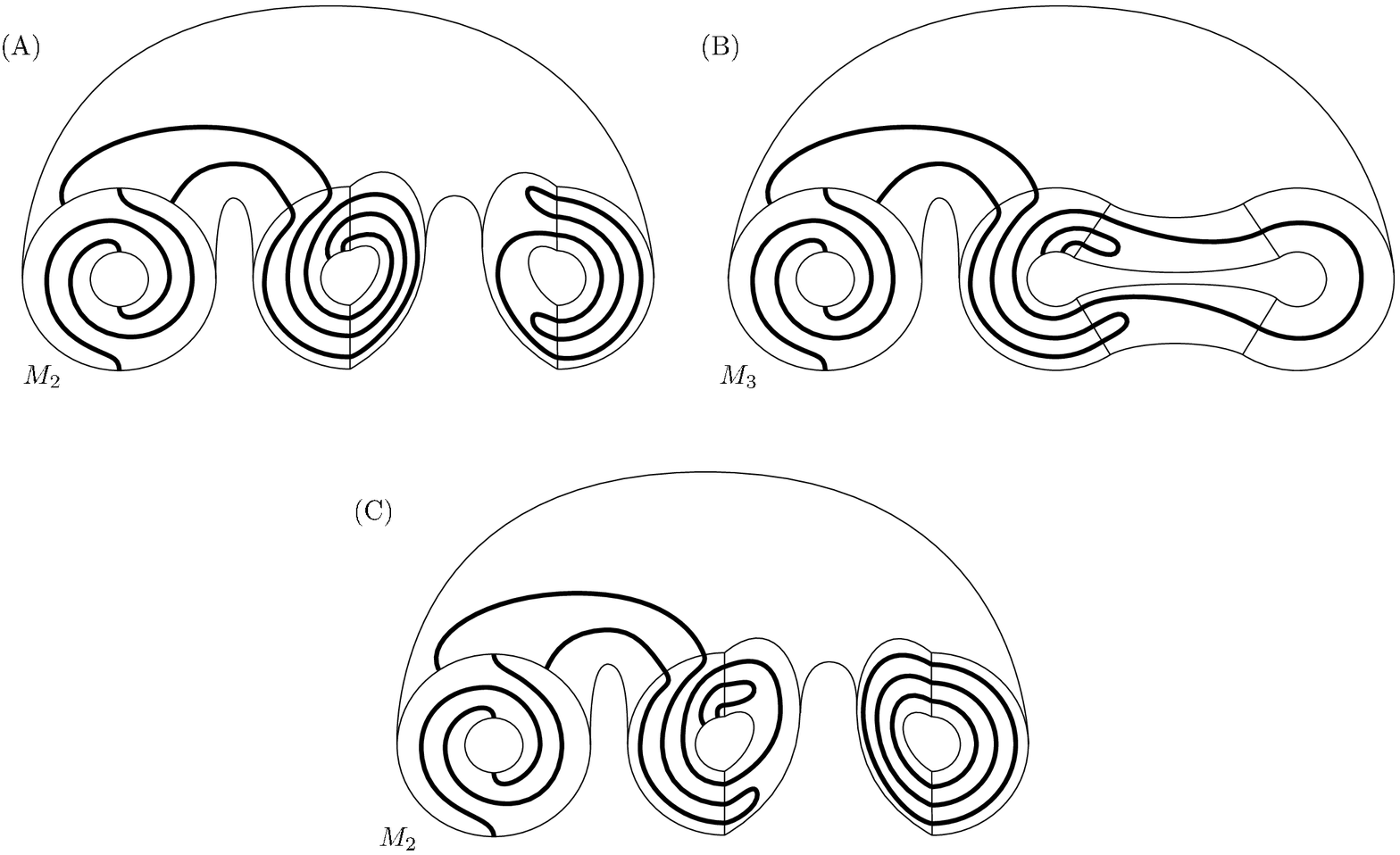}}}
	\caption{Transition from ${\mbox{\it I}}_k$ to $\ii^\pm_n$, continued.} 	
	\label{fig12}
\end{figure}

The case of $\beta_1$ occupies Figures~\ref{fig12}(A,B) and the case 
of $\beta_2$ occupies Figure~\ref{fig12}(C).  We will explain $\beta_2$ first.  
After rounding the edges, $\Gamma_{\bdry M_2\setminus A^+}$ has one 
$\bdry$-parallel arc along each boundary component of $A^+$.  This implies 
that there exists a Legendrian divide $\alpha_1$ on $\bdry M_2\setminus A^+$ 
parallel to either boundary component of $A^+$ (after possibly perturbing 
$\bdry M_2\setminus A^+$ as in LeRP).  Let $\alpha_2$ be an efficient 
Legendrian curve on $A^+$, isotopic to the core curve and 
satisfying $|\alpha_2\cap \Gamma_{A^+}|=2$.  If we take an annulus spanning 
$\alpha_1$ to $\alpha_2$, by the Imbalance Principle we will obtain a degenerate 
bypass along $\alpha_2$.  Attaching the degenerate bypass will give the 
transition to some $\ii^\pm_n$.  

On the other hand, $\beta_1$ does not immediately give rise to a  
$\bdry$-parallel arc along $\bdry A^+$.  Therefore, we cut again, this time 
along $\delta_2\times I$, which is convex with efficient Legendrian boundary, to 
obtain $M_3$.   Write $\bdry \delta_2=q_1-q_0$.  Now, $|\bdry(\delta_2\times I) 
\cap \Gamma_{\bdry M_2}|=2k+2$, and $2k+1$ of the intersections are on 
$q_0\times I$, whereas $1$ intersection is on $q_1\times I$.   If $k>1$, then 
there will always be a bypass along $q_0\times I$, which transitions us to 
${\mbox{\it I}}_{k-1}$.  On the other hand, there is one extra case when $k=1$ --- after the 
edges are rounded (Figure~\ref{fig12}(B)), there exists a $\bdry$-parallel arc 
along $\bdry A^+$ on $\bdry M_3\setminus A^+$. Therefore,  we will always have a 
state transition to some $\ii^\pm_n$, provided ${\mbox{\it I}}_k$ does not represent an 
$I$-invariant tight contact structure. 
\end{proof}

We claim that $PD(\tilde{e}([\gamma_0,\gamma_0])) = \pm 
2\gamma_0$ or $0$.  To see this, first note that if $\beta$ is any curve with 
$|\gamma_0\cap \beta|=0$, then $\langle 
\tilde{e}([\gamma_0,\gamma_0]),\beta\times I\rangle=0$, since 
$\Gamma_{\beta\times I}$ consists solely of closed curves parallel to the core 
curve.  This takes care of $2g-1$ generators of $H_1(\Sigma;\Z)$.    Next, if 
$\beta$ satisfies $|\gamma_0\cap \beta|=1$, then $\langle 
\tilde{e}([\gamma_0,\gamma_0]),\beta\times I\rangle$ is $-2$, $0$, or $2$, 
depending on the configuration of $\bdry$-parallel dividing curves on 
$\beta\times I$.  This proves the claim.

Suppose now that $[\gamma_0,\gamma_0]$ is not $I$-invariant.  Let $\gamma\subset 
\Sigma$ be a closed oriented curve satisfying $\langle 
\gamma_0,\gamma\rangle=+1$.  By Proposition~\ref{exists-bypass}, 
$[\gamma_0,\gamma_0]$ contains a basic slice, and, by the freedom of choice,  
there exists a factorization into $[\gamma_0,\gamma]\cup \ll \gamma,\gamma_0; 
\pm(\gamma_0-\gamma)\rr$.  (Notation: when we write $[a_1,a_2]\cup[a_2,a_3]\cup 
\dots\cup [a_{k-1},a_k]$, the contact structure is layered in order.)  We 
initially have the possibilities 
$PD(\tilde{e}([\gamma_0,\gamma]))=\pm\gamma_0\pm \gamma$ from 
Theorem~\ref{basecase}.   However, by the claim in the previous paragraph, if 
$PD(\tilde{e}(\ll\gamma,\gamma_0\rr))=\gamma_0-\gamma$, then  
$PD(\tilde{e}([\gamma_0,\gamma])) = \pm\gamma_0 + \gamma$.  Therefore, there is 
a total of 4 possibilities: $$[\gamma_0,\gamma_0]= [\gamma_0,\gamma; 
\pm\gamma_0+\gamma]\cup \ll\gamma,\gamma_0;\gamma_0-\gamma\rr,$$
or
$$[\gamma_0,\gamma_0]=[\gamma_0,\gamma; 
\pm\gamma_0-\gamma]\cup \ll\gamma,\gamma_0;-(\gamma_0-\gamma)\rr.$$

\begin{lemma}
All four contact structures are tight, distinct, and can be embedded in a basic 
slice. \end{lemma}

\begin{proof}
Let $\gamma_0$, $\gamma_1$ be nonseparating curves satisfying $\langle 
\gamma_1,\gamma_0\rangle=+1$.   We claim that $\ll 
\gamma_0,\gamma_1;\gamma_1-\gamma_0   \rr$ can be factored into
$[\gamma_0,\gamma_0]\cup \ll\gamma_0,\gamma_1\rr$, where the first factor is 
{\em not} $I$-invariant.    This is a consequence of Proposition~\ref{ubiquity} 
as follows.  First we layer $[\gamma_0,\gamma_1] = [\gamma_0,\gamma]\cup 
[\gamma,\gamma_1]$, where $|\gamma\cap \gamma_i|=1$, $i=0,1$, using 
Proposition~\ref{ubiquity}.  Now, applying Proposition \ref{ubiquity} to 
$[\gamma,\gamma_1]$, we expand: $$ [\gamma_0,\gamma_1]=   [\gamma_0,\gamma]\cup 
[\gamma,\gamma_0]\cup \ll\gamma_0,\gamma_1\rr.$$   The union of the first and 
second slices on the right-hand side of the equality cannot be $I$-invariant by 
the semi-local Thurston-Bennequin inequality (see \cite{Gi00b}).  Therefore, we 
have obtained a factorization $$\ll \gamma_0,\gamma_1; 
\gamma_1-\gamma_0\rr = [\gamma_0,\gamma_0] \cup \ll\gamma_0,\gamma_1\rr .$$
It remains to compute $PD(\tilde{e})$ of the second factor.  If we reconcile 
$PD=\pm 2\gamma_0$ or $0$ for the first factor with 
$PD=\pm(\gamma_1-\gamma_0)$ for the second factor, we easily see that:
$$\ll \gamma_0,\gamma_1; 
\gamma_1-\gamma_0\rr = [\gamma_0,\gamma_0; 0] \cup 
\ll\gamma_0,\gamma_1;\gamma_1-\gamma_0\rr .$$
We therefore have realized at least one tight contact structure 
$[\gamma_0,\gamma_0;0]$ which is not $I$-invariant.  We also obtain another 
non-$I$-invariant tight contact structure $[\gamma_0,\gamma_0;0]$ by starting 
from $\ll\gamma_0,\gamma_1;\gamma_0-\gamma_1\rr$ instead.

Our next claim is that the two non-$I$-invariant $[\gamma_0,\gamma_0;0]$ are 
distinct.  Suppose we further factor:
$$ \ll\gamma_0,\gamma_1;\gamma_1-\gamma_0\rr=   [\gamma_0,\gamma_1]\cup 
\ll \gamma_1,\gamma_0 \rr\cup \ll \gamma_0,\gamma_1; \gamma_1-\gamma_0\rr.$$  
The relative Euler classes on the right-hand side of the equation, in order, are 
$\pm\gamma_0\pm\gamma_1$, $\pm(\gamma_0+\gamma_1)$, $\gamma_1-\gamma_0$.  (The 
reason we have $\pm(\gamma_0+\gamma_1)$ for the second term is due to the 
relative orientations of $\gamma_0$ and $\gamma_1$.)  For the union of the 
second and third layers to be tight, the second layer must have relative Euler 
class $\gamma_0+\gamma_1$ in order to cancel the $\gamma_0$'s (and the first 
layers must be $-(\gamma_0+\gamma_1)$).  Therefore, 
$$ \ll\gamma_0,\gamma_1;\gamma_1-\gamma_0\rr=   
[\gamma_0,\gamma_1;-(\gamma_0+\gamma_1)]\cup \ll \gamma_1,\gamma_0; 
\gamma_0+\gamma_1 \rr\cup \ll \gamma_0,\gamma_1; \gamma_1-\gamma_0\rr.$$  
Applying the same calculation to $\ll\gamma_0,\gamma_1;\gamma_0-\gamma_1\rr$, we 
see that the two non-$I$-invariant tight contact structures 
$[\gamma_0,\gamma_0;0]$ can be distinguished by the factorization into 
$[\gamma_0,\gamma_1]\cup \ll \gamma_1,\gamma_0\rr$.

It remains to dig further to obtain the remaining two tight contact structures 
$[\gamma_0,\gamma_0]$ with $PD(\tilde e)=\pm \gamma_0$.   It suffices to factor
$\ll\gamma_0,\gamma_1; \gamma_1-\gamma_0\rr$ into  
$$[\gamma_0,\gamma_0;-2\gamma_0]\cup\ll \gamma_0,\gamma_1; 
\gamma_0-\gamma_1\rr\cup \ll\gamma_1,\gamma_0; \gamma_0+\gamma_1\rr\cup 
\ll\gamma_0,\gamma_1; \gamma_1-\gamma_0\rr.$$   Similarly, we obtain 
$[\gamma_0,\gamma_0; 2\gamma_0]$. \end{proof}

\subsection{Proof of Parts 1 and 2 of Theorem~\ref{classification} and of 
Theorem~\ref{gluing}}    

We first need the following lemma:

\begin{lemma}  \label{exists-basic}
If $\gamma_0\not=\gamma_1$, then $[\gamma_0,\gamma_1]$ contains a basic slice.
\end{lemma}

\begin{proof}
First suppose that $|\gamma_0\cap \gamma_1|\not=0$.  Then we can realize 
$\gamma_1\times\{0,1\}$ on $\Sigma_0\cup \Sigma_1$ by efficient Legendrian 
curves, apply the Imbalance Principle to the convex annulus with boundary 
$\gamma_1\times\{0,1\}$, and find a bypass along $\gamma_1\times\{0\}\subset 
\Sigma_0$ from the interior of $\Sigma\times I$. Let $\alpha$ be the Legendrian  
arc of attachment for the bypass.  There are two possibilities for $\alpha$: (i)
$\alpha$ starts on a dividing curve $\gamma_0^1$ (one of the two curves
parallel to $\gamma_0$ on $\Sigma_0$), passes through the other curve
$\gamma_0^2$, and ends on $\gamma_0^1$; (ii) $\alpha$ starts on $\gamma_0^1$,
passes through $\gamma_0^2$, and ends on $\gamma_0^2$ after going through a
nontrivial loop.   (i) is clearly an attachment which gives rise to a basic 
slice.  For (ii), let $P\subset \Sigma_0$ be a pair-of-pants neighborhood of the
union of $\alpha$ and the annulus bounded by $\gamma_0^1$ and $\gamma_0^2$. One 
of the boundary components is a curve $\gamma$ parallel to $\gamma_0^i$, the 
second boundary component is $B(\gamma)$, parallel to the dividing curves on 
$\Sigma_{1/2}$ obtained by isotoping $\Sigma_0$ through the bypass attached  
along $\alpha$ and the third curve denoted $\delta$ may be thought of as the 
nontrivial loop $\alpha$ goes around (see Figure~\ref{typeC}).

We claim that $\gamma=\gamma_0$ and $B(\gamma)=\gamma_{1\over 2}$ are not 
isotopic.  If they were, then they would cobound an annulus $B$, and $B\cup P$ 
would be a once-punctured torus.  In a once-punctured torus, an efficient, 
nontrivial arc or closed curve will intersect another only in positive 
intersections or only in negative intersections.  This contradicts the 
efficiency of the original curve $\gamma_1\times\{0\}$. 

Therefore, we may shrink $\Sigma \times [0,1]$ and assume that 
$\gamma_0\not=\gamma_1$ are disjoint and nonisotopic. In such a situation, let 
$\beta$ denote a curve that is efficient, intersects $\gamma_0$ and does not
intersect $\gamma_1$. By using the Imbalance Principle as above, we can again 
find a bypass of either type (i), in which case we have found a basic slice, or 
a bypass of type (ii). In the latter case we can shrink again and have 
$\Sigma\times[0,1]$ with $\Gamma_{\Sigma_0}$ parallel to $\gamma$ and 
$\Gamma_{\Sigma_1}$ parallel to $B(\gamma)$, where $\gamma$, $B(\gamma)$ and 
$\delta$ form a boundary of a pair of pants $P \subset \Sigma_0$ and $\gamma$ 
and $B(\gamma)$ are not isotopic.

Assume $\delta$ is nonseparating. If $\gamma$ and $\delta$ lie on the 
same connected component of $\Sigma_0 \setminus  int(P)$, let $\beta_1$ be an 
arc in $P$ connecting $\gamma$ and $\delta$ and let $\beta_2$ be an arc in 
$\Sigma_0\setminus int(P)$ connecting the same points on $\gamma$ and $\delta$ 
as $\beta_1$. Let $\beta = \beta_1 \cup \beta_2$ (Figure~\ref{typeCcases}, Case 
$C_1$). Since $\beta$ does not intersect $B(\gamma)$, the Imbalance Principle 
applied to it produces a bypass of type (i), and hence a basic slice.

If $B(\gamma)$ and $\delta$ lie on the same connected component of $\Sigma_0 
\setminus int(P)$, an analogous argument produces $\beta$ that intersects 
$B(\gamma)$ once and does not intersect $\gamma$, and another application of the 
Imbalance Principle produces a basic slice. (See Figure~\ref{typeCcases}, Case 
$C_2$.)

If $\delta$ is separating, denote the component of $\Sigma_0 \setminus P$ it 
bounds by $S_1$. Let $\beta_1$ be a nonseparating arc in $S_1$ which starts and 
ends on $\delta$, and let $\beta_2$ be a nonseparating arc in another component 
of $\Sigma_0\setminus P$ which starts and ends on $\gamma$. Let $\beta$ 
be a closed nonseparating curve obtained by joining those arcs by arcs in $P$ 
(Figure~\ref{typeCcases}, Cases $C_3$ and $C_4$). Note that a nonseparating  
$\beta_2$ exists and $\beta$ can be chosen efficient and Legendrian in both 
$\Sigma_i$ because $\gamma$ and $B(\gamma)$ are not isotopic.  Applying the 
Imbalance Principle to this $\beta$ produces a bypass of type (ii) with 
nonseparating $\delta$. This in turn, we have shown, contains a basic slice. 
\end{proof}

\n
{\it Note.} A similar but slightly more involved argument will be carried out
in Proposition~\ref{2parallelcurves}.  The figures we refer to are the same as 
the ones we need for that argument.

\s\n
Now that we know that $[\gamma_0,\gamma_1]$ contains a basic slice, we may apply 
Proposition~\ref{fact0}, together with Proposition~\ref{ubiquity}, to  factor: 
\begin{eqnarray*}
[\gamma_0,\gamma_1] &= &\ll\alpha_0=\gamma_0,\alpha_1; 
\alpha_1-\alpha_0\rr \cup\ll\alpha_1,\alpha_2\rr \cup \dots 
\\ & & \dots \cup \ll \alpha_{k-2},\alpha_{k-1}\rr \cup [\alpha_{k-1},\alpha_k=\gamma_1],
\end{eqnarray*}
subject to the following:

\be
\item $\alpha_i$, $i=1,\dots,k$, are oriented nonseparating curves.
\item $|\alpha_i\cap \alpha_{i+1}|=1$ ($i=0,\dots,k-1$) and 
$|\alpha_i\cap\alpha_{i+2}|=0$ ($i=0,\dots, k-2$).
\item $\langle\alpha_{i+1},\alpha_{i}\rangle =1$,  $i=0,\dots, k-1$.  
\item All the slices except for the last are basic slices.
\item Without loss of generality, $PD(\tilde 
e)$ of the first factor is $\alpha_1-\alpha_0$. 
\ee

\n
We will inductively prove that the rest of the $PD(\tilde e)$'s must be, in 
order, $\alpha_2-\alpha_1, \alpha_3-\alpha_2, ... , 
\alpha_{k-1}-\alpha_{k-2}, \pm\alpha_k-\alpha_{k-1}$.

\begin{lemma}  \label{consistency}
Suppose $[\alpha_{i-1},\alpha_{i+1}]$ is the union 
$\ll\alpha_{i-1},\alpha_i\rr\cup\ll\alpha_i,\alpha_{i+1}\rr$ of two 
basic slices with $\langle \alpha_i,\alpha_{i-1}\rangle = \langle 
\alpha_{i+1},\alpha_i\rangle=+1$, $|\alpha_{i-1}\cap\alpha_{i+1}|=0$, and 
$PD(\tilde e(\ll\alpha_{i-1},\alpha_i\rr))=\alpha_i-\alpha_{i-1}$.  If 
$[\alpha_{i-1},\alpha_{i+1}]$ is tight, then $PD(\tilde 
e(\ll\alpha_i,\alpha_{i+1}\rr))=\alpha_{i+1}-\alpha_i$.  \end{lemma}

\begin{proof}  Since $\ll\alpha_i,\alpha_{i+1}\rr$ is a basic slice, we know 
that $PD(\tilde e(\ll\alpha_i,\alpha_{i+1}\rr))=\pm(\alpha_{i+1}-\alpha_i)$.   
If $PD(\tilde e(\ll\alpha_i,\alpha_{i+1}\rr))=-\alpha_{i+1}+\alpha_i$, then 
\begin{equation}  \label{eq2}
PD(\tilde e([\alpha_{i-1},\alpha_{i+1}]))=-\alpha_{i-1}+2\alpha_{i}-\alpha_{i+1}.
\end{equation}

To obtain a contradiction, we use the fact that 
$|\alpha_{i-1}\cap \alpha_{i+1}|=0$ and calculate the possible $PD(\tilde 
e([\alpha_{i-1},\alpha_{i+1}]))$ using a different method.     If $\gamma$ is 
any closed curve with $|\gamma\cap \alpha_{i-1}|=|\gamma\cap \alpha_{i+1}|=0$, 
then $\langle \tilde e([\alpha_{i-1},\alpha_{i+1}]),\gamma\times I\rangle=0$.  
There are now two possibilities: either $\alpha_{i-1}$ and $\alpha_{i+1}$ are 
homologous or they are not.  If they are homologous, then there are $2g-1$  
generators $\gamma$ for $H_1(\Sigma;\Z)$ satisfying $|\gamma\cap 
\alpha_{i-1}|=|\gamma\cap \alpha_{i+1}|=0$ and which therefore evaluate to zero. 
If $\gamma$ is a closed curve satisfying $|\gamma\cap \alpha_{i-1}|=|\gamma\cap 
\alpha_{i+1}|=1$, then  $\langle \tilde e([\alpha_{i-1},\alpha_{i+1}]), 
\gamma\times I\rangle=\pm 2$ or $0$, depending on the signs of the 
$\bdry$-parallel components.  Hence, 
\begin{equation}  \label{eq3}
PD(\tilde e([\alpha_{i-1},\alpha_{i+1}]))=\pm \alpha_{i-1} \pm\alpha_{i+1}=\pm 
2\alpha_{i-1} \mbox{ or } 0.
\end{equation}   
Next, if $\alpha_{i-1}$ and $\alpha_{i+1}$ are not homologous, then there are 
$2g-2$ generators $\gamma$ for $H_1(\Sigma;\Z)$ satisfying $|\gamma\cap 
\alpha_{i-1}|=|\gamma\cap \alpha_{i+1}|=0$ and which therefore evaluate to zero.
There are two other basis elements $\gamma$, $\gamma'$ of $H_1(\Sigma;\Z)$ which 
satisfy $|\gamma\cap \alpha_{i-1}|=0$,   $|\gamma\cap \alpha_{i+1}|=1$, and 
$|\gamma'\cap \alpha_{i-1}|=1$,   $|\gamma'\cap \alpha_{i+1}|=0$.  We evaluate
$\langle \tilde e([\alpha_{i-1},\alpha_{i+1}]),\gamma\times I\rangle=\pm 1$ and
$\langle \tilde e([\alpha_{i-1},\alpha_{i+1}]),\gamma'\times I\rangle=\pm 1$, 
which give 
\begin{equation}  \label{eq4}
PD(\tilde e([\alpha_{i-1},\alpha_{i+1}]))=\pm \alpha_{i-1} \pm\alpha_{i+1}.   
\end{equation}            
It now suffices to note that Equation~\ref{eq2} is in contradiction with 
Equations~\ref{eq3} or \ref{eq4} --- simply intersect with $\alpha_{i-1}$.   
Therefore, we are left with $PD(\tilde 
e(\ll\alpha_i,\alpha_{i+1}\rr))=\alpha_{i+1}-\alpha_i$. \end{proof}

Thus, by Lemma~\ref{consistency}, we find that the $PD(\tilde e)$'s of the basic 
slices are $\alpha_2-\alpha_1, \alpha_3-\alpha_2, ... , 
\alpha_{k-1}-\alpha_{k-2}.$  Finally, although the last slice is not a basic 
slice, an argument almost identical to that of Lemma~\ref{consistency} proves 
that the relative Euler class is $\pm\alpha_k-\alpha_{k-1}$.  Therefore, we see 
that the initial basic slice $\ll\alpha_0,\alpha_1\rr$ uniquely                               
determines all the subsequent basic slices and reduces the possibilities for 
the last slice to two.   Thus, there are at most $4$ possibilities for 
$[\gamma_0,\gamma_1]$, up to isotopy rel boundary.  
Adding up the relative Euler classes of the slices, we obtain $PD(\tilde 
e([\gamma_0,\gamma_1]))=\pm\gamma_0\pm\gamma_1$.   The relative Euler classes 
distinguish the 4 possibilities, provided $\gamma_0$ and $\gamma_1$ are not 
homologous. 

We now have the following proposition:

\begin{prop} \label{prop:euler}
Suppose $[\alpha,\alpha']$ is a tight contact structure, where $\alpha$ and 
$\alpha'$ are nonseparating.  Then $PD(\tilde e([\alpha,\alpha']))=\pm\alpha\pm 
\alpha'$.  Moreover, if $[\alpha,\alpha']$ admits a factorization 
$\ll\alpha_0=\alpha,\alpha_1\rr\cup\ll\alpha_1,\alpha_2\rr\cup \dots \cup 
\ll\alpha_{k-1},\alpha_k=\alpha'\rr$ with $\langle \alpha_{i+1}, 
\alpha_{i}\rangle=+1$, then $PD(\tilde 
e([\alpha,\alpha']))=\pm(\alpha'-\alpha).$ \end{prop}

\begin{proof} 
This was largely proved in the above paragraphs, with the difference that we 
required that $|\alpha_i\cap \alpha_{i+2}|=0$.  This extra condition is not 
required in Proposition~\ref{prop:euler}, since there always exists a 
subdivision which satisfies this extra property. \end{proof}

Theorem~\ref{gluing} immediately follows from Proposition~\ref{prop:euler}.  The 
following two lemmas complete the proof of Theorem~\ref{classification}.

\begin{lemma} \label{findthem}
All four $[\gamma_0,\gamma_1]$ are $\Sigma\times I$ layers inside 
some basic slice $\ll\gamma_0,\gamma\rr$ with $\langle 
\gamma,\gamma_0\rangle=+1$. 
\end{lemma}   

\begin{proof}
We will start with $\ll\gamma_0,\gamma; 
\gamma-\gamma_0\rr$ and find two of the four possibilities;  the other two can 
be found inside $\ll\gamma_0,\gamma; -\gamma+\gamma_0\rr$.  Using 
Proposition~\ref{ubiquity}, we find a factorization: $$\ll\gamma_0,\gamma; 
\gamma - \gamma_0\rr=\ll\alpha_0=\gamma_0,\alpha_1\rr \cup 
\ll\alpha_1,\alpha_2\rr \cup \dots \cup 
\ll\alpha_{k-1},\alpha_k=\gamma_1\rr\cup[\alpha_k,\gamma].$$ The tight contact 
structure $[\gamma_0,\gamma_1]\subset[\gamma_0,\gamma]$, obtained by layering 
$\ll\alpha_0,\alpha_1\rr,\dots,\ll\alpha_{k-1},\alpha_k\rr$, must have 
$PD(\tilde e)=\gamma_1-\gamma_0$ by Proposition~\ref{prop:euler}.    Now, by 
Theorem~\ref{gluing}, if $PD(\tilde e([\gamma_0,\gamma]))=\gamma-\gamma_0$, then 
$PD(\tilde e([\gamma_0,\gamma_1]))$ must be $\gamma_1-\gamma_0$.  To obtain 
$-(\gamma_0+\gamma_1)$, we start with $\ll\gamma_0,\gamma;\gamma-\gamma_0\rr$ 
and factor: 
$$\ll\gamma_0,\gamma;\gamma-\gamma_0\rr=\ll\alpha_0=\gamma_0,\alpha_1\rr 
\cup \ll\alpha_1, \alpha_2\rr \cup \dots \cup 
\ll\alpha_{k-1},\alpha_k=\gamma_1\rr\cup\ll\alpha_k,\gamma\rr\cup 
\ll\gamma,\alpha_k\rr\cup [\alpha_k ,\gamma], $$ 
with all but the last layer basic.  Throwing away the last layer on the 
right-hand side, we obtain $[\gamma_0,\gamma_1]$ with $PD(\tilde 
e)=-(\gamma_0+\gamma_1)$.  
\end{proof}

\begin{lemma}[Unique factorization] \label{unique-factorization}
The 4 tight contact structures on $[\gamma_0,\gamma_1]$ are distinct.
\end{lemma}

\begin{proof}
If $\gamma_0$ and $\gamma_1$ are not homologous, then the four tight contact 
structures are distinguished by the relative Euler class.  If 
$\gamma_0\not=\gamma_1$ are homologous, then the relative Euler class cannot 
distinguish between $\gamma_1-\gamma_0$ and $\gamma_0-\gamma_1$.   We already 
showed that one of the $[\gamma_0,\gamma_1]$ admits a factorization:
$$[\gamma_0,\gamma_1; \gamma_1-\gamma_0] = \ll\alpha_0=\gamma_0,\alpha_1; 
\alpha_1-\alpha_0 \rr \cup\ll\alpha_1,\alpha_2;\alpha_2-\alpha_1\rr \cup  \dots 
\cup \ll\alpha_{k-1},\alpha_k=\gamma_1;\alpha_k-\alpha_{k-1}\rr,$$
with $\langle \alpha_{i+1},\alpha_i\rangle = +1$.  We claim that given any other 
factorization:
$$[\gamma_0,\gamma_1; \gamma_1-\gamma_0] = \ll\alpha_0=\gamma_0,\alpha_1 \rr 
\cup\ll\alpha_1,\alpha_2 \rr \cup  \dots \cup  \ll\alpha_{k-1},\alpha_k = 
\gamma_1\rr,$$  
each $\ll\alpha_i,\alpha_{i+1}\rr$ has relative Euler class 
$\alpha_{i+1}-\alpha_i$.  From the discussion in Lemma~\ref{consistency}, we 
see that it suffices to prove that the relative Euler class of 
$\ll\alpha_{k-1},\alpha_k\rr$ is $\alpha_k-\alpha_{k-1}$.   But now, by 
Lemma~\ref{findthem}, we see that if $\alpha_{k+1}$ satisfies 
$\langle\alpha_{k+1},\alpha_k\rangle=1$, then  
$[\gamma_0,\gamma_1]\cup\ll\alpha_k,\alpha_{k+1};\alpha_{k+1}-\alpha_k\rr$ is 
tight.  Now, applying Lemma~\ref{consistency}, we see that 
$\ll\alpha_{k-1},\alpha_k\rr$ must have relative Euler class 
$\alpha_k-\alpha_{k-1}$.  
\end{proof}

\section{Classification of tight contact structures on hyperbolic 3-manifolds 
which fiber over the circle}

In this section we provide the proof of Theorem~\ref{thm:fibration}. Let $M$ be 
a closed, oriented 3-manifold which fibers over the circle with fibers (oriented 
surfaces $\Sigma$) of genus $g\geq 2$ and pseudo-Anosov monodromy 
$f:\Sigma\rightarrow \Sigma$.  In other words, $M=\Sigma\times [0,1]/\sim$, 
where $(x,0)\sim (f(x),1)$.  Recall that the pseudo-Anosov condition is 
equivalent to saying that for every multicurve $\Gamma\subset \Sigma$, 
$f(\Gamma)\not=\Gamma$.  The assumption that $e(\xi)$ is extremal guarantees 
that the dividing set on any convex fiber $\Sigma$ is a union of pairs of 
parallel curves bounding annuli.  To prove Theorem~\ref{thm:fibration}, we first 
show that there exists a convex fiber $\Sigma$ for which $\Gamma_\Sigma$ 
consists of exactly two nonseparating curves.   This is accomplished by starting 
with an arbitrary convex fiber $\Sigma$ and inductively reducing the number of 
curves in $\Gamma_\Sigma$ by two, by isotoping $\Sigma$ through an appropriate 
bypass.  The following proposition will be used to show the existence of an 
appropriate bypass.

\begin{prop}\label{goodgamma} 
Let $\xi$ be a tight contact structure on $\Sigma\times [0,1]$ with convex 
boundary, and suppose $\Gamma_{\Sigma_0}\not=\Gamma_{\Sigma_1}$. Then there is a 
closed efficient curve $\gamma$, possibly separating, such that $|\gamma \cap 
\Gamma_{\Sigma_0}| \neq | \gamma \cap \Gamma_{\Sigma_1}|$. \end{prop}

\begin{proof}
Suppose the dividing set $\Gamma_{\Sigma_0}$ is the disjoint union, over 
$i=1,\dots, k$, of $m_i$ curves isotopic to $\delta_i$. 
We will then identify the dividing set $\Gamma_{\Sigma_0}$ with the point $\{m_i 
\delta_i\}_{i=1}^k$ in the weighted curve complex on $\Sigma_0$. After an 
isotopy, we may assume that all $\delta_i$ and $f(\delta_j)$ intersect 
transversely and efficiently (realize the geometric intersection number).  If 
$\delta_i \cap f(\delta_j) \neq \emptyset$ for some $i$ and $j$, then letting 
$\gamma = \delta_i$ satisfies the conclusion of the proposition.

We may now assume that $\delta_i \cap f(\delta_j) = \emptyset$ for all pairs 
$i,j=1,\dots, k$.  Thus $\{\delta_i\} \cup \{f(\delta_i)\}$ may be completed to 
a pair-of-pants decomposition of $\Sigma$.  It is not hard to show directly that 
weights on the cuffs of a pair-of-pants decomposition are determined by their 
intersection number with transverse embedded curves; thus the required curve 
$\gamma$ exists. (For more details, see \cite{H2}.)
\end{proof}

The following grew out of discussions with John Etnyre: 
 
\begin{prop}\label{2parallelcurves} 
Let $\xi$ be a tight contact structure on $M$ with $\langle e(\xi),\Sigma\rangle 
= \pm(2g-2)$.  Then there exists a convex surface isotopic to the fiber whose 
dividing set consists of two parallel nonseparating curves.
\end{prop}

\begin{proof}
Let $\Sigma_0$ be a convex surface isotopic to a fiber. Cut $M$ along $\Sigma_0$ 
and denote the cut-open manifold $\Sigma \times I $ and the contact structure 
$\xi|_{\Sigma\times I}$ by $\xi$. Then $\Gamma_{\Sigma_1} = 
f(\Gamma_{\Sigma_0})$.

By Proposition~\ref{goodgamma}, we may choose $\gamma$ to be an efficient curve 
such that there is an imbalance in the number of intersections of $\bdry(\gamma 
\times I)$ with $\Gamma_{\Sigma_0}$ and $\Gamma_{\Sigma_1}$.  (Note that if 
$\gamma$ is separating, then we may need to use the stronger form of LeRP: 
Lemma~\ref{lemma:super-lerp}, or its proof.) There is a bypass 
contained in $\Sigma\times I$, attached along either $\Sigma_0$ or $\Sigma_1$. 
We can assume without loss of generality that the bypass is along $\Sigma_0$. 
Since $\gamma$ is chosen to be efficient, the bypass can be neither trivial nor 
increase the number of dividing curves on $\Sigma_0$.

Denote the attaching arc of the bypass by $\alpha$. We will discuss the possible 
types of bypass attachments, and in each case show that the number of dividing 
curves can eventually be decreased by two.      

\s\n 
{\bf Type A.} The attaching arc $\alpha$ intersects three different dividing 
curves. 

\s\n
In this case, there exists a consecutive parallel pair of dividing curves which 
intersect $\alpha$.   By isotoping $\Gamma_{\Sigma_0}$  through 
the bypass, we will therefore reduce the number of dividing curves by two.  See 
Figure~\ref{typeA}.

\begin{figure}[ht] 		
{\epsfysize=2in\centerline{\epsfbox{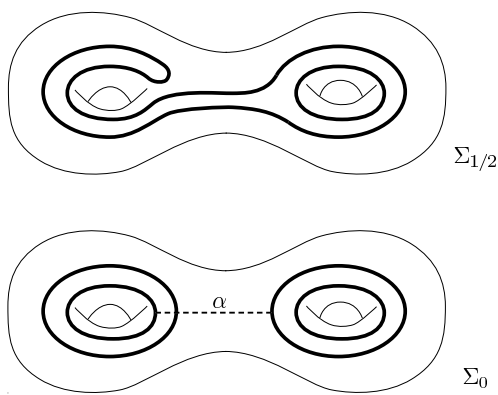}}} 	
\caption{Type A attaching arc} 		
\label{typeA} 
\end{figure}

\s\n
{\bf Type B.} The attaching arc $\alpha$ starts on a dividing curve $\gamma_1$, 
passes through a parallel dividing curve $\gamma_2$, and ends on $\gamma_1$.  

\s\n
Isotop $\Sigma_0$ through this bypass to obtain $\Sigma_{1/2}$.  Let $N \subset
\Sigma_0$ be a punctured torus regular neighborhood of the union of $\alpha$ and 
the annulus bounded by $\gamma_1$ and $\gamma_2$.  Identify the region between 
$\Sigma_0$ and $\Sigma_{1/2}$ with $\Sigma \times [0, {1\over 2}]$ in such a way 
that the contact structure is $I$-invariant on $(\Sigma \setminus N) \times [0, 
{1\over 2}]$ and is a basic slice on $N \times [0, {1\over 2}]$. See 
Figure~\ref{typeB}(A).

\vskip.1in

\begin{figure}[ht]	
	{\epsfysize=2in\centerline{\epsfbox{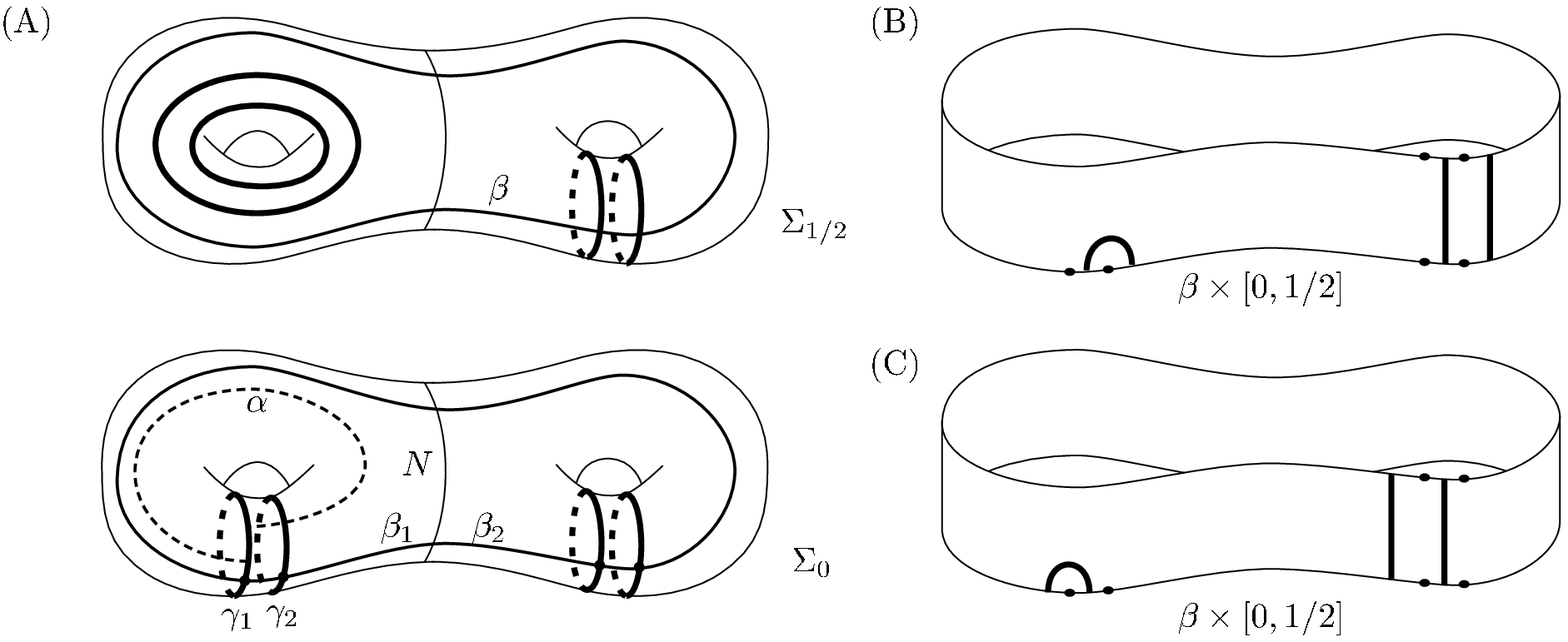}}}
	\caption{Type B attaching arc} 	
	\label{typeB}
\end{figure}

Since there are more than two dividing curves on $\Sigma_0$, there are dividing 
curves contained in $\Sigma_0 \setminus N$. Choose $\beta \subset \Sigma$ to be 
a closed curve formed out of an arc $\beta_1 \subset N$ and an arc $\beta_2 
\subset \Sigma_0 \setminus N$, where the arcs have the following properties: (i)  
$\beta_1 \times \{ 0 \}$ intersects each of $\gamma_0$ and $\gamma_1$ once, (ii)  
$\beta_1 \times \{ {1\over 2} \}$ intersects no dividing curves, (iii) $\beta_2 
\times \{ 0 \}$ intersects dividing curves in $\Sigma \setminus N$, and (iv) 
$\beta \times \{ 0 \}$ is efficient with respect to $\Gamma_{\Sigma_0}$.  We 
may choose $\beta \times [0, {1\over 2}]$ to be convex.   Since the contact 
structure is $I$-invariant on $(\Sigma \setminus N)\times [0,{1\over 2}]$, the 
dividing curves along $\beta_2$ are vertical.  Thus there are only two possible 
dividing curve configurations, both of which are shown in 
Figure~\ref{typeB}(B,C), and either of these forces the existence of a bypass of 
Type~A along a subarc of $\beta \times \{ 0\}$.

\s\n
{\bf Type C.} The attaching arc $\alpha$ starts on a dividing curve $\gamma_1$, 
passes through a parallel dividing curve $\gamma_2$, and ends on $\gamma_2$ 
after going around a nontrivial loop.

\s\n
Let $P \subset \Sigma_0$ be a pair-of-pants regular neighborhood of the union of 
$\alpha$ and the annulus bounded by $\gamma_1$ and $\gamma_2$.  One of the 
boundary components of $P$ is a curve $\gamma$ parallel to $\gamma_1$ and 
$\gamma_2$, the second boundary component is a curve $B(\gamma)$ parallel to the 
two new dividing curves obtained from $\gamma_1$ and $\gamma_2$ after isotoping 
through the bypass along $\alpha$, and the third curve $\delta$ may be thought 
of as the nontrivial loop $\alpha$ goes around. See Figure~\ref{typeC}. 

\begin{figure}[ht]	
	{\epsfysize=4in\centerline{\epsfbox{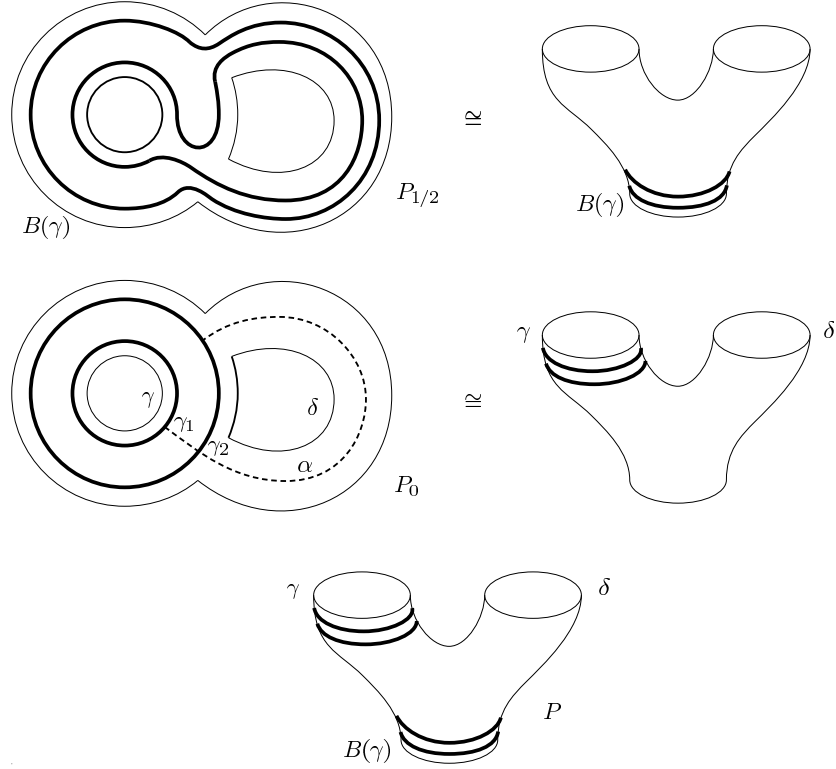}}}
	\caption{Type C attaching arc} 	
	\label{typeC}
\end{figure}

Isotop $\Sigma_0$ through this bypass to obtain $\Sigma_{{1/ 2}}$. Identify the 
region between $\Sigma_0$ and $\Sigma_{{1/ 2}}$ with $\Sigma \times [0, 
{1\over 2}]$ in such a way that the contact structure is $I$-invariant on $(\Sigma 
\setminus P) \times [0, {1\over 2}]$.

The proof now proceeds roughly as in the Type~B case, that is, a curve $\beta$ to which the
Imbalance Principle can be applied will be produced.  To do this, we must 
consider the possible ways that $P$ can sit in $\Sigma$.  See 
Figure~\ref{typeCcases}.

\begin{figure}[ht]	
	{\epsfysize=5in\centerline{\epsfbox{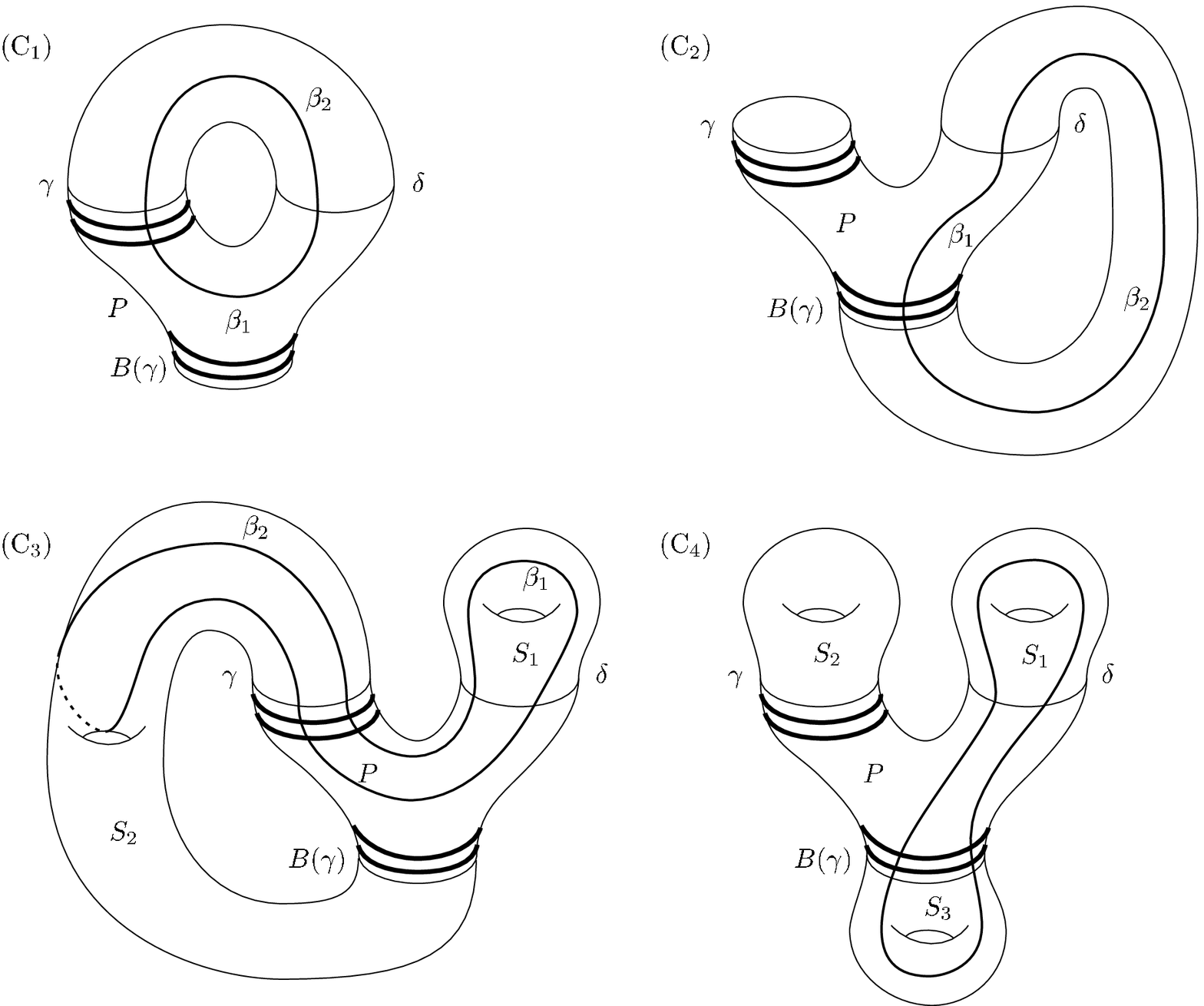}}}
	\caption{Type C cases} 	
	\label{typeCcases}
\end{figure}

\s\n                                              
{\bf Type $\mbox{C}_{\bf 1}$.} $\delta$ is nonseparating and there is an arc 
$\beta_2 \subset \Sigma \setminus P$ connecting $\gamma$ and $\delta$.

\s\n
Let $\beta_1 \subset P$ be an arc connecting the same points of $\gamma$ and 
$\delta$ and let $\beta = \beta_1 \cup \beta_2$.  The arcs may be chosen so that 
$\beta$ is efficient.  If $\beta_2$ intersects any dividing curves, then the 
Imbalance Principle applied to $\beta$ (or more precisely to $\beta \times [0, 
{1\over 2}]$) produces a bypass of Type~A.  Otherwise, since $\beta$ is 
nonseparating, $\beta \times \{{1\over 2}\}$ is nonisolating and can be made a 
Legendrian divide.  Again the Imbalance Principle applied to $\beta \times [0, 
{1\over 2}]$ produces a bypass, this time a degenerate one which can be 
perturbed to a bypass of Type~B.  

\s\n
{\bf Type $\mbox{C}_{\bf 2}$.} $\delta$ is nonseparating and there is an arc 
$\beta_2 \subset \Sigma \setminus P$ connecting $B(\gamma)$ and $\delta$.

\s\n
Let $\beta_1 \subset P$ be an arc connecting the same points of $B(\gamma)$ and
$\delta$ and let $\beta = \beta_1 \cup \beta_2$. Just as in Type~$\mbox{C}_ 
{1}$, the Imbalance Principle applied to $\beta$ produces a bypass of Type~A or 
B.    

\s\n
{\bf Type $\mbox{C}_{\bf 3}$.} $\delta$ is separating and there is an arc 
contained in $\Sigma \setminus P$ connecting $B(\gamma)$ and $\gamma$.

\s\n
Let $S_1$ be the component of $\Sigma \setminus P$ bounded by $\delta$, and let 
$S_2$ be the component with which contains $\gamma$.  The genus of $S_1$ is 
greater than 0 because $\delta$ was assumed to be nontrivial.  Recall $\alpha$ 
is a subarc of an efficient curve obtained by Proposition~\ref{goodgamma}.  If 
$S_2$ is an annulus, then the subsurface $S_2\cup P$ is a once-punctured torus 
and $\alpha$ cannot be a subarc of an efficient arc on $S_2\cup P$.  Therefore,   
$S_2$ must have genus greater than 0 also.

Since there are more than two dividing curves on $\Sigma_0$, there must be 
dividing curves contained in $S_1 \cup S_2$. Let $\beta_1 \subset S_1$ be a 
nonseparating arc starting and ending on $\delta$, let $\beta_2 \subset S_2$ be 
a nonseparating arc starting and ending on $\gamma$, and choose the arcs so that 
at least one of them has nontrivial, essential intersection with the dividing 
curve set.  Pick two disjoint arcs in $P$ connecting the endpoints of $\beta_1$ 
to the endpoints of $\beta_2$ and let $\beta$ be the union of all four arcs.  
The Imbalance Principle produces either a bypass of Type~A or, since the 
$\beta_i$ were chosen to be nonseparating, a bypass of Type~$\mbox{C}_{1}$ or 
$\mbox{C}_{2}$.

\s\n
{\bf Type $\mbox{C}_{\bf 4}$.} All three curves $\delta$, $\gamma$, and 
$B(\gamma)$ are separating.

\s\n
Let $S_1, S_2$ and $S_3$ be, in order, the components of $\Sigma \setminus P$ 
these curves bound. Since each $S_i$ must have genus greater than 0, the proof 
is the same as in Type~$\mbox{C}_{3}$, with one possible extra case.  If all of 
the dividing curves of $\Sigma \setminus P$ are in $S_3$, then, to force $\beta 
\times \{{1\over 2}\}$ to be isolating, the desired $\beta$ must be produced in 
$S_3 \cup P \cup S_1$, and the resulting bypass will lead to a reduction  in the 
number of dividing curves on $\Sigma_{1/2}$ instead of $\Sigma_0$.  

Thus in all possible cases a bypass (or a sequence of bypasses) can be found 
that will reduce the number of dividing curves by two.  Hence, we can assume 
$\Gamma_{\Sigma_0}$ consists of two parallel curves. Moreover, we can assume 
they are nonseparating. If they are not, we can find a nonseparating  curve 
intersecting $\Gamma_{\Sigma_0}$  and not intersecting $\Gamma_{\Sigma_1}$ and 
use  the inbalance principle to find a bypass that transforms 
$\Gamma_{\Sigma_0}$  into a pair of nonseparating curves. \end{proof}

\begin{proof}[Proof of Theorem~\ref{thm:fibration}]
Fix a nonseparating curve $\gamma$ on $\Sigma$.  Since $f$ is pseudo-Anosov, 
$\gamma$ and $f(\gamma)$ are not isotopic, and we orient these curves so that 
$f$ preserves their orientations.

Let $\xi$ be the universally tight contact structure on $\Sigma \times [0,1]$ 
with $PD(\tilde e(\xi)) = f(\gamma) - \gamma$.  We claim that the contact 
structure obtained by gluing $\Sigma\times\{0\}$ with $\Sigma\times\{1\}$ via 
$f$ is universally tight.   It is immediate, by the Gluing Theorem~\ref{gluing}, 
that $2n$ copies of $\Sigma \times [0,1]$ stacked and glued together via $f$ 
will produce a universally tight contact structure on $\Sigma \times [-n,n]$.  
It follows that the glued-up contact structure on $M$ is universally tight, for 
any potential overtwisted disk in $M$ would lift to the $\Sigma \times \R$ cover 
of $M$ and therefore would be contained in some $\Sigma \times [-n,n]$.

To prove uniqueness, use Proposition~\ref{2parallelcurves} to choose a convex 
fiber $\Sigma'$ such that $\Gamma_{\Sigma'} = 2\gamma'$, where $\gamma'$ is some 
nonseparating curve, and split $M$ along $\Sigma'$ to obtain $\Sigma' \times I$.  
Since $f$ is pseudo-Anosov, $\Gamma_{\Sigma' \times 
\{0\}}\not = \Gamma_{\Sigma'\times\{1\}}$, and therefore $\Sigma'\times I$ 
contains a basic slice by Lemma~\ref{exists-basic}.  It follows from 
Proposition~\ref{ubiquity} that a new convex fiber $\Sigma''$ may be chosen such 
that $\Gamma_{\Sigma''} = 2\gamma$.

Let us assume for simplicity that $\gamma$ and $f(\gamma)$ are not 
homologous.  (The general case is no harder, but requires the use of 
Lemma~\ref{unique-factorization}, and is harder to state.) By 
Theorem~\ref{classification}, there are four possibilities for $\xi_{\Sigma'' 
\times I}$ with distinct relative Euler classes $\pm f(\gamma) \pm \gamma$.  If 
the relative Euler class is $\pm f(\gamma) + \gamma$, then $\Sigma'' \times 
[0,1]$ contains a slice $\Sigma''\times [0,{1\over 2}]$ where $\Gamma_{\Sigma'' 
\times \{{1\over 2}\}}=2\gamma$ and the slice has relative Euler class 
$2\gamma$.  By peeling this slice and reattaching using $f$, we obtain 
$\Sigma''\times[{1\over 2},1]\cup f(\Sigma''\times [0,{1\over 2}])$.  It 
follows that for any given tight contact structure on $M$, there always exists a 
convex fiber $\Sigma$ such that the restriction of $\xi$ to $\Sigma \times I$ 
(obtained by cutting along $\Sigma$) has relative Euler class $\pm f(\gamma) - 
\gamma$.

It is clear, from the location of bypasses on the $-f(\gamma) - \gamma$ contact 
structure on $\Sigma \times I$ (or by Proposition~\ref{consistency}), that 
gluing by $f$ produces an overtwisted contact structure on $M$.  Thus $\xi$ on 
$M$ can always be split to have relative Euler number $f(\gamma) - \gamma$ and 
the uniqueness on $M$ follows from  Theorem~\ref{classification}.

It remains to discuss weak fillability.  First note that there exists a 
$C^0$-small perturbation of the fibration into a universally tight contact 
structure $\xi$, by a result of Eliashberg and Thurston \cite{ET}.  Since the 
Euler class evaluated on the fiber is unchanged under the perturbation, we have, 
say,  $\langle e(\xi),\Sigma\rangle=-(2g-2)$.  By the uniqueness which we just 
proved, the unique extremal tight contact structure on $M$ with $\langle 
e,\Sigma\rangle=-(2g-2)$ is isotopic to $\xi$.  Now, to prove weak symplectic 
fillability, we construct a symplectic 4-manifold $(X,\omega)$ with $\bdry X=M$ 
for which $\omega|_{T\Sigma}>0$ (for all fibers $\Sigma$).  Since $\xi$ is close 
to the fibration, we would then be done.  The construction of $X$ is relatively 
straightforward from the Lefschetz fibration perspective, once we observe that 
any element $f:\Sigma\rightarrow \Sigma$ of the mapping class group of a closed 
surface $\Sigma$ can be written as a product of positive Dehn twists. (For more 
details on symplectic Lefschetz fibrations, see, for example \cite{GS}.) We take 
a symplectic Lefschetz fibration $X\rightarrow D^2$ with generic fiber $\Sigma$ 
and a singular fiber for each positive Dehn twist in the product expression.  We 
then see that $\bdry X=M$ has the desired monodromy and each fiber $\Sigma$ is a 
symplectic submanifold.  \end{proof}

\s\n
{\it Acknowledgements.}   The first author gratefully acknowledges the American 
Institute of Mathematics and IHES for their hospitality.  The authors also thank 
Emmanuel Giroux for suggesting that Theorem~\ref{thm:fibration} might be true, 
and John Etnyre for discussions which led to Proposition~\ref{2parallelcurves}.

\end{document}